\documentclass[11pt,amssymb]{article}
\usepackage{amsfonts,latexsym}

\newcommand{\RR}{\mathbb R} 
\newcommand{\NN}{\mathbb N}

\newcommand{\ZZ}{\mathbb Z}
\newcommand{\R}{\mathbb R} 
 
\newcommand{\D}{\mathcal D}
\newcommand{\Q}{\mathcal Q}

\newcommand{\del}{\partial} 
\newcommand{\e}{\varepsilon} 
 
\newcommand{\ha}{\frac{1}{2}}
\newcommand{\al}{\alpha}
\newcommand{\tB}{\tilde{B}}
\newcommand{\bC}{\mbox{\bf C}}
\newcommand{\Se}{\Sigma_\e}
\newcommand{\dSe}{\del \Sigma_\e}

\newcommand{\tSeP}{\tilde{\Sigma}_{\e,{\cal P}}} 
\newcommand{\tLeP}{\tilde{L}_{\e,{\cal P}}}
\newcommand{\x}{\mbox{\bf x}}

\newcommand{\Le}{L_{\e}} 
 
\newcommand{\Lej}{{L}_{\e,j}}

\newcommand{\tw}{\tilde{w}}

\newcommand{\bw}{\bar{w}} 
\newcommand{\bu}{\bar{u}}

\newtheorem{theorem}{Theorem} 
\newtheorem{lemma}{Lemma} 
\newtheorem{proposition}{Proposition} 
\newtheorem{corollary}{Corollary} 
\newtheorem{definition}{Definition} 
\newtheorem{remark}{Remark} 
 
\begin{document} 
 
\title{Constant mean curvature surfaces \\ with Delaunay ends} 
 
\author{Rafe Mazzeo \thanks{Supported by 
the NSF under Grant DMS-9626382 and a Young Investigator Fellowship, 
and by a Sloan Foundation Fellowship} \\ Stanford University 
\and Frank Pacard \thanks{Supported by 
the American Institute of  Mathematics} \\ Universit\'e Paris XII } 

\maketitle 
 
\section{Introduction and statement of the results} 
 
In this paper we shall present a construction of Alexandrov-embedded 
complete surfaces $M$ in ${\RR}^3$ with finitely many ends and finite
topology, and with nonzero constant mean curvature (CMC). This construction
is parallel to the well-known original construction by Kapouleas \cite{K1}, but
 we feel that ours somewhat simpler analytically, and controls the resulting 
geometry more closely. On the other hand, the surfaces we construct have 
a rather different, and usually simpler, geometry than those of Kapouleas; 
in particular, all of the surfaces constructed here are noncompact, so we do 
not obtain any of his immersed compact examples. The method we use here 
closely parallels the one we developed recently \cite{MP1} to study the 
very closely related problem of constructing Yamabe metrics on 
the sphere with $k$ isolated singular points, just as Kapouleas' construction 
parallels the earlier construction of singular Yamabe metrics by Schoen 
\cite{S1}.  

The original examples of noncompact CMC surfaces were those in the 
one-parameter family of rotationally invariant surfaces discovered by 
Delaunay in 1841 \cite{D}.
One extreme element of this family is the cylinder; the `Delaunay surfaces' 
are periodic, and the embedded members of this family interpolate between 
the cylinder and an infinite string of spheres arranged along a common axis. 
The family continues beyond this, but the elements now are immersed, and we 
shall not consider them here. We are mostly interested in surfaces which 
are {\it Alexandrov embedded}. This condition, by definition, means that 
although the surface may be immersed, the immersion extends to one of the 
solid handlebody bounded by the surface. The CMC surfaces we construct
here have Alexandrov embedded ends; if the minimal $k$-noids out of which 
they are built (as we describe below) are Alexandrov embedded, then the
whole surfaces satisfy this condition.

The r\^ole of  Delaunay surfaces in the theory of complete CMC surfaces is 
analogous to the r\^ole of catenoids (and planes) 
in the study of complete minimal surfaces of finite total curvature. 
For example, just as any complete 
minimal surface with two ends must be a catenoid \cite{S2},
it was proved by Meeks \cite{Me} and Korevaar, Kusner  and Solomon \cite{KKS} 
that any constant mean curvature surface with at most two ends is necessarily 
a Delaunay surface.  A rather more remarkable theorem, paralleling the fact 
that any end of a complete minimal surface of finite total curvature must 
be asymptotic to a catenoid or a plane, 
is the fact that any end of a CMC surface of the type we consider must be  
asymptotic to one of these rotationally symmetric Delaunay surfaces (and in
particular, must be cylindrically bounded). 

The fact that such CMC surfaces exist in abundance was proved, as noted 
above, by Kapouleas \cite{K1} in 1987. More recently, using much softer 
methods based on the Schwarz reflection principle, Grosse-Brauckmann has 
constructed families of CMC surfaces with $k$ ends and with $k$-fold 
dihedral symmetry \cite{KGB}. 

The general analysis of the moduli space of CMC surfaces was considered by 
the first author, Kusner and Pollack \cite{KMP} (essentially merely 
translating the analogous results in \cite{MPU} for the singular Yamabe 
problem). The basic result is that the moduli
space ${\cal M}_{g,k}$ of $k$-ended surfaces with genus $g$ is a locally real
analytic variety of virtual dimension $3k$ (before dividing out by the action
of the group of Euclidean motions). This virtual dimension is attained at
any point $\Sigma \in {\cal M}_{g,k}$ where a certain analytic nondegeneracy
criterion is satisfied. This condition will be explained in some detail
below. It seems very difficult to  decide whether any of the surfaces 
constructed by Kapouleas satisfy this nondegeneracy criterion. 
This was one motivation for the present work, because the solutions we 
construct do satisfy it. We note finally that recently, Kusner, 
Grosse-Brauckmann and Sullivan gave a heuristic determination of the moduli 
space ${\cal M}_{0,3}$; this work strongly suggests that all elements of 
this space are nondegenerate. In work in progress, they expect to
make their arguments rigorous.

We now state our main result in more detail. This result is simply that CMC 
surfaces may be constructed out of certain building blocks in a specified 
and controlled manner. 
There are only two types of building blocks: Delaunay surfaces (or more 
precisely, halves of Delaunay surfaces) and minimal $k$-noids. The former 
we have already encountered; on the other hand, a minimal $k$-noid is by 
definition a complete minimal surface $\Sigma$ of finite total curvature 
and with $k$ ends. Denote the moduli space of minimal $k$-noids of genus 
$g$ by ${\cal H}_{g,k}$. This space has been studied by Perez and 
Ros \cite{PR}, and they prove the result corresponding to that of \cite{KMP} 
that this space is real analytic of virtual dimension $3k$. 
Elements of various of these spaces have been shown to exist by the classical
Weierstrass method, and more recently elements have been constructed for $g$
very large by Kapouleas by his gluing methods \cite{K2}.  Again there is
a notion of nondegeneracy of such surfaces, and a surface $\Sigma$ is a smooth
point in its corresponding moduli space precisely when it satisfies this
nondegeneracy condition. Fortunately, the Weierstrass representation gives
sufficient information in some cases to establish the existence of 
nondegenerate minimal $k$-noids. Quite recently, Cosin and Ros \cite{CR}
have proven the existence of nondegenerate minimal $k$-noids with
specified weight parameters (in the sense described later here).
Finally, as noted earlier, any end of an 
element $\Sigma \in {\cal H}_{g,k}$ is asymptotic to the end of a plane 
or a catenoid. We define ${\cal H}^0_{g,k}$ to be the subset of 
${\cal H}_{g,k}$ consisting of nondegenerate
elements $\Sigma$ with all ends of $\Sigma$ catenoidal, rather than planar. 

We may now state the main 
\begin{theorem} 
Fix any $\Sigma_0 \in {\cal H}^0_{g,k}$. Then there exists a family of CMC 
surfaces $\Sigma_\e \in {\cal M}_{g,k}$ constructed by gluing half-Delaunay 
surfaces onto
each end of the dilated surface $\e\Sigma_0$. The surfaces $\Sigma_\e$ have
the property that for any $R > 0$, the dilated surface $\e^{-1}\Sigma_\e$
restricted to the ball $B_R(0)$ converges in the ${\cal C}^\infty$ topology 
to the restriction of $\Sigma_0$ to $B_R(0)$. In addition, $\Sigma_\e$ are 
regular points of ${\cal  M}_{g,k}$.
\label{th:1.1} 
\end{theorem} 

Our proof has some novel features.  Rather than finding solutions
as perturbations off of degenerating families of approximate solutions,
as has been common in such constructions, we instead find (infinite
dimensional) families of CMC surfaces as normal graphs over each of
the component pieces, the half-Delaunay surfaces and a truncated $k$-noid. 
By studying the Cauchy data of the functions producing these normal graphs, 
and prove that we may match the Cauchy data from the inner piece
(the graph over the $k$-noid) with that from the ends, thus one advantage
of this procedure is that more of the technical complications caused
by the nonlinearities are avoided than would be otherwise. 

The plan of this paper is as follows. We first discuss the Delaunay surfaces in
some detail, collecting and proving various technical properties concerning
them that we require later, specifically those concerned with their behavior
in the singular limit, as they approach the bead of spheres. This is followed 
by the analysis of the Jacobi (i.e. linearized CMC) operator, especially in 
this singular limit, for the half-Delaunay surfaces. We then use this to 
discuss the full family of CMC surfaces in a neighborhood of these 
rotationally invariant surfaces, as usual, keeping careful track of the 
behavior in the limit. We then turn to a discussion of minimal surfaces with 
$k$ catenoidal ends, i.e. the $k$-noids, briefly reviewing their geometry and 
then  treating the relevant aspects of the linear analysis of their Jacobi 
operator.
After that we can approach the family of CMC surfaces obtained as normal
graphs over suitable truncations of these $k$-noids. At last we can
put all of this together and prove that it is possible to match the
Cauchy data, and so obtain the proof of the main theorem.

\noindent
{\bf Aknowledgment~:} This paper was written when the second author was 
visiting the University of Stanford. The second author would like to take 
here the opportunity to thank the American Institute of Mathematics and 
the Mathematics Departement for their support and hospitality.
 
\section{Notation, conventions and definitions} 

We recall some basic facts about the geometry of immersed surfaces, 
and review various ways the equations for constant mean curvature
may be specified. Some good references for this material are the
book by Osserman \cite{O} and the survey article by Wente \cite{W}.

Suppose that $\Sigma$ is given as the image of a regular immersion 
$\x: {\cal U} \longrightarrow {\RR}^3$. Here ${\cal U}$ is an open set in 
${\RR}^2$ with coordinates $u = (u_1,u_2)$. The unit normal to $\Sigma$ is 
defined to be
\[
\nu(u_1,u_2) = \frac{\del_{u_1}{\x} \wedge \del_{u_2}{\x}}{\|\del_{u_1}{\x} 
\wedge  \del_{u_2}{\x}\|}, 
\]
and the components of the first and second fundamental forms $g$ and $B$ 
are then
\begin{eqnarray}
E = \langle \del_{u_1}\x  , \del_{u_1}\x \rangle, \ F = \langle 
\del_{u_1}\x , \del_{u_2} \x \rangle, \ G = \langle \del_{u_2}\x, \del_{u_2}
\x \rangle \nonumber \\[2mm]
L = \langle \del^2_{u_1 u_1} \x, \nu \rangle, \ M = \langle \del^2_{u_1 u_2}
\x, \nu \rangle, \ N = \langle \del^2_{u_2 u_2}\x, \nu \rangle.
\label{eq:2.1}
\end{eqnarray}
The principal curvatures $k_1$ and $k_2$ are the eigenvalues of $B$ relative
to $g$. The mean curvature is defined to be the sum ({\it not} the average) of
the principal curvatures, $H \equiv k_1 + k_2$, and the Gauss curvature
is their product, $K \equiv k_1k_2$. 

We shall almost always be using orthogonal parameterizations (that is 
to say, parameterizations for which $F =0$), in which case the 
formul\ae\ for $H$ and $K$ reduce to:
\[
H=\frac{L+N}{E}, \qquad   K=\frac{LN-M^2}{E^2}.
\]

The important equations of surface theory are the Gauss and Codazzi
equations, which link the intrinsic and extrinsic geometry of $\Sigma$. 
Rather than write these down in general, we consider the special case
where the parameterization is isothermal, so that $E=G \equiv e^{2\omega}$
 and $F = 0$.
Let $u = u_1 + iu_2$ be the corresponding complex coordinate, and define the 
Hopf differential
\begin{equation}
\phi(u) \, du^2 = \left( \ha(L - N) - iM\right)du^2.
\label{eq:2.2}
\end{equation}
The principal curvatures are then given by
\[
k_1 = \frac{H}{2} -|\phi| e^{-2w}, \qquad k_2 = \frac{H}{2} + |\phi| e^{-2w}.
\]
If $\Sigma$ is CMC, so $H$ is constant, then the Codazzi equations are
equivalent to the holomorphy of this differential. The Gauss equation is 
simply
\begin{equation}
\Delta \omega + \frac{H^2}{4} e^{2\omega} - |\phi|^2 e^{-2\omega} = 0.
\label{eq:2.3}
\end{equation}

\section{Delaunay surfaces} 
 
We now make a detailed study of the first of the basic building blocks
we use later, the Delaunay surfaces of revolution. 

\subsection{Definition and basic equations}

The Delaunay surfaces mentioned in the introduction are surfaces of
revolution, and so we use cylindrical coordinates. In particular, if the 
axis of rotation is the vertical one, and if $t$ is a linear coordinate
along this axis and $\theta$ is the angular variable around it, then
we consider surfaces $\Sigma$ given by the parametrization
\begin{equation}
\x(t,\theta) = (\rho(t)\cos \theta, \rho(t)\sin \theta , t).
\label{eq:2.4}
\end{equation}
The condition that such a surface has constant mean curvature $1$ 
gives an ordinary differential equation for the function $\rho(t)$,
and solutions of this equation correspond to the Delaunay surfaces.

To obtain this ODE, first note that the unit normal of $\Sigma$ at $\x(t,\theta)$ is
\[
\nu(t,\theta) = \frac{1}{\sqrt{1+\rho_t^2}}(-\cos \theta,-\sin\theta, \rho_t),
\]
where subscripts denote derivatives, and then that the metric tensor and
second fundamental form are given by
\begin{equation}
g = (1 + \rho_t^2)\,dt^2 + \rho^2 \,d\theta^2, \qquad
B =- \frac{\rho_{tt}}{\sqrt{1+\rho_t^2}}\, dt^2 + 
\frac{\rho}{\sqrt{1+\rho_t^2}}\,d\theta^2.
\label{eq:3.1}
\end{equation}
It follows that the mean curvature is given by the expression
\begin{equation}
H = -\rho_{tt}(1+\rho_t^2)^{-3/2} + \rho^{-1}(1+\rho_t^2)^{-1/2},
\label{eq:5prime}
\end{equation}
and so the condition $H = 1$ becomes the equation
\begin{equation}
\rho_{tt} - \frac{1}{\rho}(1+\rho_t^2) + (1 + \rho_t^2)^{3/2} = 0.
\label{eq:3.2}
\end{equation}

There are two special solutions of (\ref{eq:3.2}) that can be
determined immediately. The first is the constant solution
$\rho_1(t) \equiv 1$, the cylindrical graph of which is the
cylinder of radius $1$. The other, $\rho_0(t) = \sqrt{4-(t-2)^2}$,
$|t-2| \leq 2$, corresponds to the sphere of radius $2$ centered at $(0,0,2)$. 
The singular limit of the Delaunay surfaces mentioned in the introduction 
corresponds to the periodic extension of $\rho_0(t)$ to all of $\RR$. 

For $0 < \e < 1$, we {\it define} $\rho_\e(t)$ to be the solution of 
(\ref{eq:3.2}) which attains its minimum value $\rho_\e(0) = \e$ at $t=0$.
By differentiating, we see that if $\rho$ is a solution of (\ref{eq:3.2})
then 
\[
{\cal H}(\rho, \rho_t) \equiv \rho^2 -\frac{2 \rho}{\sqrt{1+\rho_t^2}},
\]
is constant. In particular, ${\cal H}(\rho_\e,(\rho_\e)_t) = \e \, (\e-2) < 0$. 
Introduce the new parameter $\tau$ by $\tau^2 = \e \,(2-\e)$, so that 
$\e = 1-\sqrt{1-\tau^2}$ and $0<\tau<1$ as well. We then deduce
immediately the 
\begin{proposition} The solution $\rho_\e$ of (\ref{eq:3.2}) with
$\rho_\e(0) = \e$, $\del_t \,\rho_\e (0) = 0$ is periodic and varies between the 
limits 
\[
\e \equiv 1 - \sqrt{1-\tau^2} \le \rho_\e \le 1+ \sqrt{1-\tau^2}.
\]
In particular, $\rho_\e(t) \le 2$ for all $t$ and $\e$. 
\label{pr:3.01}
\end{proposition}

These solutions and their translates constitute the (embedded) Delaunay 
family; the surfaces determined by them, as well as their images under 
Euclidean motions, are the Delaunay surfaces. 

To simplify notation, we often drop the subscript $\e$ (which should 
not be confused with 
the standard partial derivative notation). We also
introduce a new parameterization, changing both the independent and 
dependent variables, which simplifies the study of the $\rho_\e$.
A change of independent variable corresponds to the introduction
of a function $t = k(s)$, which should be a diffeomorphism of 
$\RR$ onto itself. the function $k$ is chosen so that the corresponding 
parameterization
\[
(s,\theta) \longrightarrow (\rho(k(s))\cos\theta, \rho(k(s))\sin\theta,k(s)),
\]
is isothermal. This corresponds to the condition
\begin{equation}
(1 + \rho_t^2)\, k_s^2 = \rho^2.
\label{eq:3.3}
\end{equation}
With the initial condition $k(0) = 0$ and noting that $k_s >0$, then $\rho$
uniquely determines $k$. Also, $k_s \neq 0$ (so long as  $\rho \neq 0$, 
which is always the case here), and so using the periodicity  
of $\rho$ we see that $k$ must be a diffeomorphism. Now, use the parameter 
$\tau \in (0,1)$ from above and define the function $\sigma(s)$ by 
\begin{equation}
\rho(k(s)) = \tau e^{\sigma(s)}. 
\label{eq:3.4}
\end{equation}
A brief calculation shows that 
\begin{equation}
\rho_t = \frac{\sigma_s}{\sqrt{1-\sigma_s^2}}, \quad
1+\rho_t^2 = \frac{1}{1-\sigma_s^2}, \qquad \mbox{and}\quad
\rho_{tt} = \frac{\sigma_{ss}}{\tau e^{\sigma}(1-\sigma_s)^2}.
\label{eq:3.5}
\end{equation}
In terms of the new variable $s$ and function $\sigma$, the first
and second fundamental forms are now
\begin{equation}
g = \tau^2e^{2\sigma}\left(ds^2 + d\theta^2\right), \qquad
B = -\frac{\sigma_{ss}\tau e^{\sigma}}{\sqrt{1-\sigma_s^2}}\, ds^2
+ \tau e^{\sigma}\sqrt{1-\sigma_s^2}\,d\theta^2,
\label{eq:3.6}
\end{equation}
and (\ref{eq:3.2}) becomes
\begin{equation}
\sigma_{ss} + \tau e^{\sigma} \sqrt{1-\sigma_s^2} - (1 - \sigma_s^2) = 0.
\label{eq:3.7}
\end{equation}

We can now see that this parameterization is indeed simpler.
\begin{proposition} 
The function $\sigma$ defined by (\ref{eq:3.4}) satisfies the equation
\begin{equation}
\sigma_{ss} + \frac{\tau^2}{2} \sinh{2\sigma} = 0, 
\label{eq:3.8}
\end{equation}
and in fact 
\begin{equation}
\sigma_s^2 + \tau^2 \cosh^2 \sigma = 1.
\label{eq:3.9}
\end{equation}
Conversely, if $\sigma$ satisfies (\ref{eq:3.9}) for some 
$0 < \tau < 1$ (and hence also (\ref{eq:3.8})), 
with \newline $\tau^2 \cosh^2 \sigma (0)=1$, 
and if $t = k(s)$, where $k(0)=0$ and 
\begin{equation}
k_s(s)  = \frac{\tau^2}{2} (1 + e^{2\sigma}),
\label{eq:3.81}
\end{equation}
then $\rho(t) \equiv \tau e^{\sigma(s)}$ satisfies (\ref{eq:3.2})
and $\rho(0) = \e$ where $\e =1 - \sqrt{1-\tau^2}$.
\label{pr:3.1}
\end{proposition}
{\bf Proof:} Let $\tilde{\sigma}= \sigma +\log \tau$; since $\sigma$ solves
(\ref{eq:3.7}) then $\tilde{\sigma}$ solves
\begin{equation}
\tilde{\sigma}_{ss} + e^{\tilde{\sigma}} \sqrt{1-\tilde{\sigma}_s^2} - 
(1 - \tilde{\sigma}_s^2) = 0.
\label{eq:3.9p}
\end{equation}
The function 
\[
\tilde{\cal H}(\tilde{\sigma}, \tilde{\sigma}_s)\equiv
2 e^{\tilde{\sigma}} \sqrt{1 - \tilde{\sigma}_s^2} -e^{2\tilde{\sigma}}
\]
is a positive constant, say $\gamma$, when $\tilde{\sigma}$ is a solution of 
(\ref{eq:3.9p}). Evaluating at $s=0$, since $\sigma_s(0) =
\tilde{\sigma}_s (0) = 0$ and 
$e^{\tilde{\sigma}(0)} = \tau e^{\sigma(0)} = \e$, we have
$\gamma = 2\e - \e^2 \equiv \tau^2$. Reverting back to $\sigma$, this 
implies (\ref{eq:3.9}), and hence (\ref{eq:3.8}) by differentiation. 

The converse, that starting from $\sigma(s)$ and $\tau$, then defining 
$t = k(s)$ as in the statement of the Proposition, the corresponding
function $\rho(t)$ satisfies (\ref{eq:3.2}) is a straightforward calculation
which we leave to the reader. Notice that then Gauss and Codazzi equations are 
automatically fulfilled. \hfill $\Box$

\begin{remark}
Translating back to the notation of \S 2, the log of the conformal factor $\omega$ 
and the norm of the coefficient function of the Hopf differential are given
by 
\[
\omega = \sigma + \log \tau \qquad \mbox{and}\qquad |\phi| =\displaystyle
\frac{\tau^2}{2}, 
\]
cf. (\ref{eq:2.2}) and (\ref{eq:2.3}).
\end{remark}
Henceforth, the functions $\rho(t)$, $\sigma(s)$ and $k(s)$ will always
be related in the manner dictated by this Proposition; furthermore, the 
dependence on the parameters $\e$ and $\tau$ will not always be written 
explicitly, but we shall use them interchangeably.
We shall call either of these parameters the {\it necksize} of the
corresponding Delaunay solution. 
 
\subsection{Uniform estimates for Delaunay solutions in the singular limit} 

In this section we present a series of technical lemmata regarding the
behavior of various quantities associated to the Delaunay solutions as $\e$ 
(or $\tau$) tends to zero. Some of the estimates below are easier to obtain 
for $\rho$ and some for $\sigma$, and we shall use these functions 
interchangeably. We first estimate the period of $\sigma$; the corresponding 
estimate for $\rho$ is not required later so we merely state it and refer to 
\cite{K1} for its proof. Then we obtain some simple `global' estimates
for $\rho$, which are rather weak, but frequently useful, as well as a 
corresponding simple estimate for $\sigma$. Finer estimates for $\rho$ when 
$t$ is not too large then lead to a good comparison between the variables $s$ 
and $t$. 

\begin{proposition} Let $S_\e$ and $T_\e$ denote the periods
of $\sigma$ and $\rho$, respectively. Then as functions of
$\tau$ and $\e$,
\[
S_\e = -4 \log \tau + O(1) = -2\log \e + O(1)
\]
and 
\[
T_\e = 4 + \tau^2\log (1/\tau) + O(\tau^2) = 4 + 2 \, \e\,
\log (1/\e) + O(\e).
\]
\label{pr:est-period}
\end{proposition}
{\bf Proof:} As stated above, we only check the statement about $S_\e$.
First, using (\ref{eq:3.9}), we see that 
\[ 
\frac{1}{4} S_\e = \int_{\sigma(0)}^0 \frac{1}{\sqrt{1 -\tau^2 \cosh^2 x}}
\,dx. 
\] 
Expand the denominator into exponentials and change variables, setting
$u = e^x$. Then, letting
\[
A(\tau) \equiv e^{\sigma(0)} = \frac{1}{\tau} - \sqrt{\frac{1}{\tau^2} - 1}
= \frac{\tau}{2} + O(\tau^3),
\]
this becomes
\[
\frac{1}{4} S_\e = \int_{A(\tau)}^1 \frac{1}{\sqrt{u^2 - 
\frac{1}{4}\tau^2(u^2 +1)^2}}\,du.
\]
Changing variables once again, reduces this to an integral of the form
\[
\int_{\tau /2}^{1 - O(\tau^2)} \frac{1}{\sqrt{v^2 - \tau^2/4}}\,dv.
\]
Finally, this last integral may be computed explicitly, and
equals $-\log \tau + O(1)$. The estimate for $S_\e$ in terms
of $\e$ follows from the relationship between $\e$ and $\tau$. \hfill $\Box$ 
\smallskip

To place the next result into context, note that one limiting
solution of the basic equation (\ref{eq:3.2}) is
$\rho_0(t) = \sqrt{4-t^2}$, the equation for the sphere.
Also, the catenoid of necksize $\e$ (which is a solution of the equation
corresponding to (\ref{eq:5prime}) when $H=0$) is given
as a cylindrical graph by the function $\rho_c(t) = \e \cosh (t/\e)$.
The function $\rho(t)$ may be compared to each of these solutions.

\begin{proposition} For any $\e \in (0,1)$, the Delaunay solution
$\rho (t)= \rho_\e(t)$ satisfies the following bounds 
\begin{equation}
\e \le \rho(t)  \leq  \e \cosh (t/\e)  \,
\label{eq:3.12}
\end{equation}
\begin{equation}
1 + \rho_t^2 \  \leq  \displaystyle{\frac{1}{\e^2}}\,  \rho ^2 
\qquad \mbox{\rm (comparison with the  equation of a catenoid)}  
\label{eq:3.13}
\end{equation}
\begin{equation}
\rho^2 (1+\rho_t^2) \  \leq  4  \qquad \mbox{\rm (comparison with the 
equation of a sphere)} 
\label{eq:3.14}
\end{equation}
for any $t \in \RR$.
\end{proposition} 
{\bf Proof:} It is clear from (\ref{eq:3.9}) that $\sigma$ is monotone 
increasing on $[0, S_\e /2]$ and monotone decreasing on $[S_\e /2 ,S_\e]$. 
Correspondingly, $\rho$ is monotone increasing
on $[0, T_\e /2]$ and monotone decreasing on $[T_\e /2 ,T_\e]$. In particular,
its value at $0$ is its absolute minimum, so the lower bound in 
(\ref{eq:3.12}) is valid. 

Now multiply (\ref{eq:3.2}) by $2\rho_t(1+\rho_t^2)^{-1}$ and integrate
to get 
\begin{equation}
\log (1+\rho_t^2(t)) = 2\log(\rho(t)/\e) - 2\int_0^{t}
\rho_t(u)\sqrt{1+\rho_t^2(u)}\,du.
\label{eq:3.145}
\end{equation}
Since $\rho_t \geq 0$ for $t \in [0, T_{\e}/2]$, we get (\ref{eq:3.13}).
Next, because $\rho \geq \e$, we may write $\rho(t) =
\e\cosh(w(t)/\e)$. Inserting this into (\ref{eq:3.13}) leads to the 
inequality $w_t^2 \leq 1$. Since $w(0)=0$, we conclude that 
$w(t)\leq t$ for all $t \in [0, T_\e /2]$ and the second part
of (\ref{eq:3.12}) follows by periodicity. 
 
Since the final estimate (\ref{eq:3.14}) is not required later, we shall
not prove it here. \hfill $\Box$ 
\smallskip
 
Now we come to the more refined estimates for $\rho$.
\begin{proposition} 
For any $\e \in (0,1)$ and $t \in \RR$, the Delaunay solution $\rho (t) = \rho_\e (t)$ satisfies 
\[ 
|\rho(t) - \e \cosh (t/\e)| \leq c\, \e^2 \, e^{3 |t|/\e}
\qquad \mbox{\rm and } \qquad |\rho_t(t)-\sinh(t/\e)| \leq  c \, \e \, e^{3|t|/\e}. 
\] 
These estimates are nontrivial only when $|t| \leq 
\displaystyle{\frac{\e}{2}\log  (\frac{c}{\e})}$ for some constant $c>0$. 
\label{pr:ref-rho}
\end{proposition} 
\noindent {\bf Proof:} From (\ref{eq:3.12}) and (\ref{eq:3.13}) we get 
\[ 
1 +\rho_t^2 (t)\leq \frac{\rho^2(t)}{\e^2}\leq \cosh^2(t/\e),
\] 
and hence $|\rho_t(t)| \leq \sinh(t/\e)$ for any $t>0$. Now use this in 
(\ref{eq:3.145}) to obtain for $t>0$
\[ 
-2 \int_0^{t} \sinh (u/\e) \cosh (u/\e) \,du \leq \log(\e^2
\frac{1+\rho_t^2}{\rho^2})(t) \leq 0,
\] 
where the last inequality follows from (\ref{eq:3.13}); equivalently, 
\[ 
\rho^2(t)\left(\exp(-\e\sinh^2 (t/\e))-1\right) \leq \e^2 (1+\rho_t^2(t))
 -\rho^2(t) \leq 0. 
\] 
As we did above, set $\rho =\e\cosh(\displaystyle{{w}/{\e}})$. Then since we
already know that $w(t)\leq t$, and since $u \longrightarrow \coth u$ is 
monotone decreasing, we conclude that
\[ 
\frac{\cosh^2 ({w(t)}/{\e} )}{\sinh^2({w(t)}/{\e})} \geq  
\frac{\cosh^2 ({t}/{\e} )}{\sinh^2({t}/{\e})}. 
\] 
The previous inequality now implies that 
\[ 
\left(\exp(-\e \sinh^2 (t/\e))-1\right)\frac{\cosh^2(t/\e)}{\sinh^2 (t/\e)}
\leq ({w}_{t}^2(t)-1) \leq 0. 
\] 
Since $e^{-x} -1 \geq -x$ for all $x \geq 0$, the left side of this last 
inequality is certainly larger than or equal to $-\e \cosh^2 (t/\e)$. Hence, 
we finally obtain 
\begin{equation} 
(1-\e \cosh^2 (t/\e))^{1/2} \leq w_t(t) \leq 1, 
\label{eq:3.16}
\end{equation} 
for all $t \in [0, - \e / 2  \, \log \e ]$. Now integrate this
inequality, using that $(1-x)^{1/2} \geq 1-x$ when $0 < x < 1$ and 
$\cosh^2 x \leq e^{2|x|}$ for all $x\in \RR$, to get
\[ 
-\frac{\e^2}{2} e^{2t/\e} \leq - \e^2 \int_0^{t/\e} \cosh^2 (s)ds 
\leq w_{\e}(t)-t \leq 0. 
\] 
We conclude that for $0 < t < - \e / 2 \,   \log \e $, we have 
\[ 
\e \cosh(t/\e -\frac{\e}{2} e^{2 t/\e}) \leq \rho(t) \leq \e \cosh (t/\e). 
\] 
The estimate for $\rho$ follows at once. To get the estimate for $\rho_t$, use 
(\ref{eq:3.16}) and the relationship $\rho_t(t) =w_t(t)
\sinh(w(t)/\e)$. \hfill $\Box$ 

We can finally give a quantitative estimate for the relationship between the 
variables $t$ and $s$, or equivalently, for the function $k(s)$. 
\begin{proposition} 
For $|s| <  S_\e / 8$, the function $t = k(s)$ admits
the expansion 
\[ 
k(s) = \e s + \frac{\e^2}{8} e^{2 s} + O( \e^2 \log \e),
\] 
uniformly as $\e \rightarrow 0$. 
\label{pr:keps} 
\end{proposition} 
{\bf Proof:} It suffices to consider the case $t = k(s) \ge 0$. The estimate 
for $\rho$ from the last Proposition implies
\[
\rho(k(s))^2 = \frac{1}{4}\e^2 e^{2k/\e} + O(\e^2) + O(\e^3 e^{4k/\e}) + 
O(\e^4 e^{6k/\e}).
\]
The errors here are all of size no greater than $O(\e^2)$ precisely
when $ |k| \leq - \e/ 4 \, \log \e $ (up to an additive constant).
Assuming that $s_0 > 0$ is chosen so that this bound is satisfied, then
by the definition of $\tau$ and $\sigma$,
\[
\tau^2 e^{2\sigma(s)} = \frac{\e^2}{4} e^{2k(s)/\e} + O(\e^2).
\]
Now recall the definition of $k$ via its derivative from
Proposition~\ref{pr:3.1}, 
\[
k_s = \frac{\tau^2}{2} + \frac{\tau^2}{2}e^{2\sigma} = \e + \frac{\e^2}{8}
e^{2k/\e} + O(\e^2).
\]
Since $k \le (\e/4)\log (1/\e)$, we obtain $e^{2k/\e} \le \e^{-1/2}$, and so
$k_s = \e + O(\e^{3/2})$ in this range. Integrate to get
$k = \e s + O(\e^{3/2}s)$ for $0 \leq s \leq s_0$. From this equation,
we see that $k \le (\e/4)\log(1/\e)$ provided $s \leq S_\e /8 $,
so that we may take $s_0$ to be this last value. Now use this formula
for $k$ in terms of $s$ in the estimate for $k_s$ above to get that
\[
k_s(s) = \e + \frac{\e^2}{8} e^{2s} + O( \e^{2}  \log \e),
\]
for $0 \leq s \leq S_\e /8$. Integrating this, at last, gives
the estimate of the Proposition. \hfill $\Box$
\medskip

Collecting the results of Proposition \ref{pr:ref-rho} and the result of 
Proposition \ref{pr:keps}, we obtain~:
\begin{proposition}
There exists a constant $c>0$ independent of $\e$ such that the following
inequalities hold
\[
\tau e^\sigma \geq c \, \e^{3/4} , \qquad \qquad
\tau^2 \cosh (2\sigma)  \leq c \, \e^{1/2} , \qquad \mbox{if} 
\qquad s \in [\frac{S_\e}{8} , \frac{3S_\e}{8} ],
\]
and 
\[
\tau^2 \cosh (2\sigma) \leq c \, \e^{2} \, e^{2s} , \qquad \mbox{if} 
\qquad s \in [\frac{3 S_\e}{8} , \frac{S_\e}{2}  ] .
\]
\label{pr:af}
\end{proposition}
{\bf Proof:} Recall that $S_\e =-2 \log \e +O(1)$. It follows from Proposition
 \ref{pr:keps} that
\[
k(s) = \e s + O(\e^{3/2}) \qquad \mbox{if} \qquad s \in [0,\frac{S_\e}{8}]. 
\]
Therefore, using Proposition \ref{pr:ref-rho}, we obtain the expansion
\begin{equation}
\tau e^{\sigma} = \e \cosh s + O(\e^{3/2}e^{s}) \qquad \mbox{and} \qquad
\tau e^{-\sigma} = \frac{2}{\cosh^2 s} + O(\e^{1/2}e^{-s}) ,
\label{eq:colexp}
\end{equation}
if $s \in [0, S_\e /8 ]$. Since $\sigma$ is increasing in $[0, S_\e/2]$, we 
conclude that
\[
\tau e^{\sigma (s)} \geq \tau e^{\sigma(S/8)} \geq c \e^{3/4},
\]
for all $s\in [S_\e/8, S_\e/2]$. Similarly, since $|\sigma|$ is increasing 
in $[0, S_\e/4]$ and decreasing in $[S_\e/4, S_\e/2]$, we get
\[
\tau^2 \cosh (2 \sigma (s)) \leq \tau^2 \cosh (2 \sigma (S_\e/8))
\leq c \e^{1/2},
\]
for all $s\in [S_\e/8, 3 S_\e /8]$. Finally, since we always have 
\[
\sigma (s) =- \sigma (S_\e/2 -s),
\]
it follows at once that, for all $s \in [S_\e/2 - S_\e/8, S_\e/2 ]$, we have 
\[
\tau^2 \cosh (2 \sigma (s)) = \tau^2 \cosh (2 \sigma (S_\e/2 -s)) \leq 
c e^{2s-S} \leq \e^2 e^{2s}.
\]
This ends the proof of the Proposition. \hfill $\Box$

We finally come to some simple estimates for $\sigma(s)$. The first is
that 
\begin{equation}
\tau^2 \cosh 2\sigma \leq 2 - \tau^2.
\label{eq:3.15}
\end{equation}
This follows trivially from multiplying $\cosh(2\sigma) = 2\cosh^2 \sigma - 1$ 
by $\tau^2$ and applying (\ref{eq:3.9}).

Next, define $\xi(s) \equiv \tau \cosh \sigma(s)$. This function is periodic 
of period $S_\e /2$, attains its maximum value $\sup \xi = 1$ at $s=0$, and 
its minimum $\inf \xi = \tau$ at $s=S_\e/4$. In addition, it is a 
solution of the equation 
\begin{equation} 
\xi_{ss} =(1 +\tau^2) \xi - 2 \xi^3,
\label{eq:3.20} 
\end{equation}
which satisfies
\begin{equation} 
\xi_s^2 =(\xi^2 -\tau^2)(1-\xi^2).
\label{eq:3.201} 
\end{equation} 

\begin{proposition} 
Suppose that $s_\ell$ is any sequence of real numbers, and  that $\tau_\ell
\rightarrow 0$. Let $\sigma_\ell$ denote the function $\sigma$ when $\tau = 
\tau_\ell$, we define $\xi_\ell(s) = \tau_{\ell}\cosh \sigma_{\ell}
(s+s_\ell)$ and  $\tilde{\xi}_\ell(s) = \tau^2_{\ell}\cosh (2 \sigma_{\ell})
(s+s_\ell)$. Then  there exists an $s_0 \in \RR$ and  subsequences of the 
$\xi_\ell$ and  $\tilde{\xi}_\ell$ which either converge uniformly to 
$0$ or else converge respectively  to $1/\cosh(s+s_0)$ and $2/\cosh^2(s+s_0)$,
 uniformly on compact sets in $\RR$. 
\label{pr:3.7} 
\end{proposition} 
{\bf Proof :} Since $(\xi_\ell)_s^2 =(\xi_\ell^2 -\tau_\ell^2)(1-\xi_\ell^2)$, 
and $|\xi_\ell| \le 1$, we see that $\xi_\ell$ is bounded in ${\cal C}^1 
(\RR)$. Using (\ref{eq:3.20})  see that $\xi_\ell$ is bounded in ${\cal C}^2 
(\RR)$. This allows us to extract a subsequence which converges uniformly on 
compact subsets of $\RR$ to a solution of 
\[ 
\xi_{ss} =(1 +\tau^2) \xi - 2 \xi^3
\] 
which satisfies
\[
\xi_s^2 = \xi^2 (1-\xi^2).
\]
For $\xi_\ell$, the claim follows since the only solutions of these equations 
are $\xi \equiv 0$ or $\xi = 1/\cosh(s+s_0)$ for some $s_0\in \RR$. 
Finally, for $\tilde{\xi}_\ell$, it is sufficient to notice that
\[
\tilde{\xi}_\ell =2 \xi_\ell^2 -\tau^2_\ell,
\] 
and the claim follows.
\hfill $\Box$ 
\medskip

Because $\xi(s)$ attains its supremum at $s=0$ we next conclude that
\begin{corollary} 
As $\e \rightarrow 0$ the families of functions $\tau\cosh\sigma(s)$ and
$\tau^2 \cosh(2\sigma(s))$ converge to $1/\cosh s$ and 
$2/\cosh^{2}s$, respectively, uniformly on compact sets. 
\label{cor:3.1} 
\end{corollary} 

In fact, we may improve the range on which the convergence in this
last Corollary takes place.
\begin{corollary}
As $\e \rightarrow 0$,
\[
\tau \cosh \sigma(s) = 1/\cosh s + O(\e^{1/2}), \qquad
\tau^2 \cosh 2\sigma(s) = 2/\cosh^2 s + O(\e^{1/2}),
\]
uniformly for $|s| \le  S_\e /8$.
\label{cor:3.2}
\end{corollary}
{\bf Proof:} Note that $\xi(s) = \del_s k(s)/\rho(k(s))$. The estimates
here follow from inserting the estimates for $k(s)$ and $\rho(t)$ 
from Propositions \ref{pr:ref-rho} and \ref{pr:keps} above.
\hfill $\Box$

Finally, since $\xi(s)$ is decreasing on $[0, S_\e/4]$ and increasing on
$[S_\e/4,S_\e/2]$, we obtain using the previous Corollary,  the
\begin{proposition} 
For all $\eta >0$, there exists an $\e_0 \in (0, 1)$ and an $s_0 >0$ such that
whenever $0 < \e \leq \e_0$ and $N/2 \,  S_\e + s_0 \leq s \leq 
(N+1)/2 \, S_\e - s_0$ for some $N \in \ZZ$, then 
\[ 
\xi(s) = \tau \cosh \sigma \leq \eta \qquad \mbox{and} \qquad \tilde{\xi} = 
\tau^2 \cosh (2 \sigma) \leq \eta.
\] 
\label{pr:3.6} 
\end{proposition} 

\section{The Jacobi operator on degenerating Delaunay surfaces}

In this section we first give an explicit expression for the linearization of 
the mean curvature operator about any one of the Delaunay surfaces, and
then proceed to develop its Fredholm theory on weighted H\"older spaces.
This theory was already developed for weighted Sobolev spaces in \cite{KMP},
and the results are essentially identical. In particular, we need to find 
spaces on which this Jacobi operator is surjective. As usual, we also 
need this surjectivity with as good control as possible as the necksize
shrinks. The results and proofs here are very close to those in \cite{MP1}.

\subsection{The Jacobi operator} Recall from the last section the 
cylindrical parameterization $\x_\e(t,\theta)$
for the Delaunay surface $\Sigma_\e$ of necksize $\e$, and the corresponding
expression for its unit normal $\nu_\e$. Given any function $w(t,\theta)$ on 
$\Sigma_\e$, its normal graph
\begin{equation}
\begin{array}{rlll}
{\x}_{w}(t, \theta ) & = & \x_\e(t,\theta) + w(t,\theta)\nu(t,\theta)\\[2mm]
\, & = & \displaystyle \left( \left(\rho - \frac{w}{\sqrt{1+\rho_t^2}}\right)
\cos \theta, \left( \rho - \frac{w}{\sqrt{1 + \rho_t^2}}\right) \sin \theta,
t + \frac{w}{\sqrt{1 + \rho_t^2}}\right) ,
\end{array}
\end{equation}
gives a regular parametrization of a surface $\Sigma_w$, provided
$w$ is sufficiently small. In terms of the coefficients of the
first and second fundamental forms of this surface, the nonlinear
operator we are interested in takes the form
\begin{equation} 
{\cal N}(w) = 1  - \frac{L_w G_w-2M_w F_w +N_w E_w}{E_w G_w-F_w^2}.
 \label{eq:4.1} 
\end{equation} 
It is well known that the linearization ${\cal L}_\e$ of ${\cal N}$ at $w=0$, 
which
is usually called the Jacobi operator for $\Sigma_\e$, is given by 
\begin{equation}
{\cal L}_\e = \Delta_{\Sigma_\e} + |A_{\Sigma_\e}|^2.
\label{eq:4.2}
\end{equation}
In terms of the parameterization above, this may be written as
\begin{equation} 
{\cal L}_\e = \frac{1}{\rho\sqrt{1+ \rho_t^2}}\del_t\left(\frac{\rho}{
\sqrt{1+\rho_t^2}}
\del_t \right)  + \frac{1}{\rho^{2}} \del^2_\theta  +  
\frac{\rho^2\rho_{tt}^2 + (1+\rho_t^2)^2}{\rho^2(1+\rho_t^2)^3}.
\label{eq:4.3} 
\end{equation} 
This looks complicated, but fortunately, becomes simpler 
in the $(s,\theta)$ coordinate system introduced above. Now
\begin{equation}
{\cal L}_\e = \frac{1}{\tau^2e^{2\sigma}} \left( \del_s^2 + \del_\theta^2
+ \tau^2 \cosh (2\sigma) \right).
\label{eq:4.4}
\end{equation}
Removing the factor $(\tau^2 e^{2\sigma})^{-1}$, it will be sufficient to 
study the operator 
\begin{equation} 
\Le w = \del_s^2 + \del^2_\theta + \tau^2 \cosh (2\sigma). 
\label{eq:4.5} 
\end{equation} 
 
Our main goal now is to study the boundary problem 
\begin{equation}
\left\{
\begin{array}{rlllll}
\Le w  & = & f \qquad &\mbox{\rm in }\quad [s_0, +\infty)\times S^1 \\[2mm]
    w  & = & \phi(\theta) \qquad &\mbox{\rm on } \quad \{s_0 \} \times S^1 ,
\end{array}  \right.
\label{eq:4.6}
\end{equation}
uniformly down to $\e = 0$. Because of the rotational invariance of the 
operator $\Le$, we may introduce the eigenfunction decomposition with
respect to the cross-sectional Laplacian $\del_\theta^2$. In this way we
obtain operators 
\[
L_{\e,j} = \del_s^2 + (\tau^2 \cosh 2\sigma - j^2), \qquad 
j \in {\mathbb Z}.
\]
Since we wish to deal only with real-valued functions, we shall use
the eigenfunctions $\chi_j(\theta) = (1/\sqrt\pi)\cos (j\theta )$ for $j > 0$, 
$\chi_j(\theta) = (1/\sqrt{\pi})\sin (j\theta )$ for $j < 0$, and 
$\chi_0(\theta) = 1/\sqrt{2\pi}$.
It will frequently be useful to separate out the operators corresponding to
the indices $j = -1, 0, 1$ from the rest, and we shall often use the notation
$\Le'$ to refer to the projection of the operator acting on these three components, 
and $\Le''$ to refer to the operator acting on all
the others together. This division is natural because, by
(\ref{eq:3.15}), the term of order zero in $L_{\e,j}$ is strictly negative
when $|j| > 1$, and so the estimates for $\Le''$ follow easily from the 
maximum principle, but this is false when $|j| \leq 1$ and $\tau$ is small. 

\subsection{Jacobi fields} A deeper reason for the separation into low
and high eigencomponents in the Jacobi operator becomes apparent when one 
examines 
the Jacobi fields, i.e. the 
solutions of $\Le \phi = 0$. Any such function may be expanded into its 
eigenseries, 
$\phi = \sum \phi_j(s) \chi_j(\theta)$, and then each $\phi_j$ solves $\Lej 
\phi_j = 0$. 
It turns out that the solutions for this problem when $j = 0, \pm 1$ may be 
determined 
explicitly in terms of the functions $\rho$ or $\sigma$; in fact, these Jacobi 
fields
correspond to quite explicit one-parameter families of CMC surfaces of which
$\Sigma_\e$ is an element. To exhibit these, first note that any smooth 
one-parameter
family $\Sigma(\eta)$ of CMC surfaces, with $\Sigma(0) = \Sigma_\e$, will
have differential at $\eta = 0$ which is a Jacobi field on $\Sigma_\e$. (This 
is meant in the sense that $\Sigma(\eta)$ should be written, for small $\eta$,
as a normal graph over $\Sigma(0)$. This is possible over any fixed compact
set of $\Sigma_\e$ for some nontrivial range of values of $\eta$ which might
diminish to zero as the compact set grows. However, this is sufficient to make 
sense 
of the derivative at $\eta = 0$.) The one-parameter families of CMC surfaces 
here are simple to describe: the first two families, corresponding to the two 
different 
solutions of $L_{\e,0}\phi = 0$, arise from varying the necksize parameter 
$\e$, and translating the $t$-variable, i.e. translating along the axis of 
$\Sigma_\e$. We denote the associated Jacobi fields by $\Psi_\e^{0,-}$ and 
$\Psi_{\e}^{0,+}$,
respectively. The other families arise from either translating or rotating
the axis of $\Sigma_\e$, so that one such translation and one such rotation
will correspond to solutions $\Psi_\e^{j,+}$ and $\Psi_\e^{j,-}$ of
$\Lej$ for $j=1$, while the translation and rotation in the orthogonal
direction corresponds to solutions for $j=-1$. In fact, we may determine
these solutions explicitly in terms of the function $\rho(t)$.
Let us write
\[
\Psi_\e^{j,\pm}(t, \theta) \equiv \Phi_\e^{j,\pm}(t)\chi_j(\theta) .
\]
We obtain :
\begin{proposition} 
The coefficient functions $\Phi_\e^{j,\pm}(t)$ of the Jacobi fields 
$\Psi_\e^{j,\pm}(t, \theta)$  for $\Sigma_\e$ for $j = -1, 0, 1$, are given 
by the formul\ae\ 
\[
\begin{array}{lllll}
\Phi_{\e}^{0,+}(t) & = & \rho_t/\sqrt{1+\rho_t^2}, \\[2mm]
\Phi_{\e}^{1,+}(t) & = & \Phi_{\e}^{-1,+}(t), \\[2mm]
                   & = &-1/\sqrt{1+\rho_t^2},
\end{array}
\qquad 
\begin{array}{llll}
\Phi_{\e}^{0,-}(t) & = & -\del_\e\rho/\sqrt{1+\rho_t^2},\\[2mm]
\Phi_{\e}^{1,-}(t) & = & \Phi_{\e}^{-1,-}(t), \\[2mm]
                   & = & -(t + \rho\rho_t)/\sqrt{1+\rho_t^2}. 
\end{array}
\]
\label{pr:4.1} 
\end{proposition} 
{\bf Proof:} First consider the families given by translations.
Suppose that the function $w$ is chosen locally so that its normal graph 
is a (small) translation of magnitude $d$ of $\Sigma_\e$ along the $x$ axis. 
Thus, for some value $t'$ near to $t$ and some value of $\theta'$ near to 
$\theta$, 
\[ 
\left\{ 
\begin{array}{rllll} 
\rho(t)\cos\theta +d  & = & \rho(t')\cos\theta' - w(t',\theta')\cos\theta'
\\[2mm] 
\rho(t)\sin\theta     & = & \rho(t')\sin\theta' - w(t',\theta')\sin\theta'
\\[2mm] 
t                     & = & t' + w(t',\theta')\rho_t(t')  ,
\end{array} 
\right. 
\] 
This system is equivalent to 
\[ 
\rho^2(t)  = \rho^2(t') - 2w(t',\theta')\rho(t') + 2 d w(t', \theta') 
\cos \theta'+ w^2 (t',\theta') -  2d\rho(t')\cos\theta' + d^2
\] 
and 
\[ 
t = t' + w(t',\theta')\rho_t(t'). 
\] 
After inserting the value of $t$ from this second equation into the first and 
collecting the lowest order terms we get
\[ 
w(t',\theta')= - d \, \cos\theta' / (1+ \rho^2_t(t')) + \mbox{higher order terms}. 
\] 
Recalling that the normal of $\Sigma_\e$ at ${\x}_\e(t',\theta')$ is 
\[ 
\nu(t',\theta') = \frac{1}{\sqrt{1+\rho_t^2}}(-\cos\theta',-\sin\theta',\rho_t(t')),
\] 
we get the stated expression for $\Phi_\e^{1,+}$. The expression for
$\Phi_\e^{-1,+}$, corresponding to translations along the $y$ axis, is
derived in an identical manner. 

In fact, nearly identical arguments work in all other cases as well. 
The relevant systems of equations are
\[ 
\left\{ 
\begin{array}{rllll} 
\rho(t)\cos\theta  & = & \rho(t')\cos\theta' - w(t',\theta')\cos\theta'\\[2mm] 
\rho(t)\sin\theta  & = & \rho(t')\sin\theta' - w(t',\theta')\sin\theta'\\[2mm] 
t + d              & = & t' + w(t',\theta')\rho_t(t'),
\end{array} 
\right. 
\]  
for the translation of size $d$ along the $z$ axis, 
\[ 
\left\{
\begin{array}{rlll} 
\rho_{\e +d}(t)\cos\theta  & = & \rho_{\e}(t')\cos\theta' - 
w(t',\theta')\cos\theta' \\[2mm] 
\rho_{\e +d}(t)\sin\theta  & = & \rho_{\e}(t')\sin\theta' - 
w(t',\theta')\sin\theta' \\[2mm] 
t                          & = & t' + w(t',\theta'){(\rho_\e)_t}(t') ,
\end{array} 
\right. 
\]
for the variation of necksize, and 
\[ 
\left\{
\begin{array}{rlll} 
(\cos d)\rho(t)\cos\theta + (\sin d)\,t  & = & \rho(t')\cos\theta' - 
w(t',\theta')  \cos\theta' \\[2mm] 
\rho(t)\sin\theta  & = & \rho(t')\sin\theta' - w(t',\theta')\sin\theta'\\[2mm] 
(\cos d)\,t - (\sin d)\rho(t)\cos\theta  & = & t' + w(t',\theta')\rho_t(t'),
\end{array} 
\right. 
\] 
for a rotation of size $d$ of the $z$-axis toward the $x$-axis, and similarly
for the rotation toward the $y$-axis. 

The calculations proceed as in the first case, and we leave the details
to the reader. \hfill $\Box$ 

\begin{corollary}
The expressions for these Jacobi fields in terms of the functions $\sigma(s)$
and $k(s)$ and the parameter $\tau$ are
\[
\Phi_\e^{0,+}(s)  =   \sigma_s , 
\]
\[
\Phi_\e^{0,-}(s)  =  \displaystyle \frac{\sqrt{1-\tau^2}}{\tau} \sigma_s \,
\del_\tau k - \sqrt{1-\tau^2}\, e^{\sigma} \, \cosh \sigma \, (1 + \tau 
\del_\tau \, \sigma),
\]
\[
\Phi_\e^{1,+}(s)  =   \Phi_\e^{-1,+}(s) = -\tau \, \cosh \sigma ,
\]
\[
\Phi_\e^{1,-}(s) =  \Phi_\e^{-1,-}(s)
=  -  k \, \tau \, ( \cosh \sigma  + \sigma_s(s) \, e^{\sigma } ).
\]
\label{cor:4.1}
\end{corollary}
The proof involves simply inserting the expressions for $\rho$, $t$ and
$\e$ in terms of $\sigma$, $k$ and $\tau$ into the previous formul\ae.

We shall require later the limits of these Jacobi fields as $\e$ tends
to zero. 
\begin{proposition}  
Let $I\subset {\RR}$ be any compact interval. Then the following
limits exist uniformly for $s \in I$:
\[ 
\lim_{\e \rightarrow 0}\Phi_{\e}^{0,+}(\e s) = \tanh s, 
\qquad
\lim_{\e \rightarrow 0}\Phi_{\e}^{0,-}(\e s) =  - (1 - s \tanh s),
 \]
\[
\lim_{\e \rightarrow 0}\Phi_{\e}^{1,+}(\e s) = \lim_{\e \rightarrow 0} 
\Phi_{\e}^{-1,+}(\e s) =  -\displaystyle{\frac{1}{\cosh s}},
\]
\[
\lim_{\e \rightarrow 0}\displaystyle{\frac{1}{\e}}
\Phi_{\e}^{1,-}(\e s) = 
\lim_{\e \rightarrow 0}\displaystyle{\frac{1}{\e}}
\Phi_{\e}^{-1,-}(\e s) =  -(\displaystyle{\frac{s}{\cosh s}} + \sinh s).
\]
\label{pr:4.2} 
\end{proposition} 
{\bf Proof:} We have used the variable $s$ in the statement of the
theorem because by Proposition~\ref{pr:keps}, $\e s/t \rightarrow 1$
uniformly for $s \in I$. The limits may be calculated using either the 
estimates 
for $\rho$ from Proposition~\ref{pr:ref-rho}, or else the expressions for 
these Jacobi fields from Corollary~\ref{cor:4.1} in terms of 
$\sigma$, and then using the limiting behaviour of $\sigma(s)$
as determined in Corollary~\ref{cor:3.1}.
\hfill $\Box$ 

The Jacobi fields we have considered so far, $\Psi_\e^{j,\pm}$, $j = 0, \pm 1$,
are all either bounded (in fact periodic) or linearly growing (because both 
$k$ and $\del_\tau k$ are linearly growing). There are 
 of course, two linearly independent solutions of the equation $\Lej \phi = 0$ 
for all $j$ with $|j|>1$ as well. It is proved in \cite{KMP}, following 
\cite{MPU},
that there exists a discrete sequence of positive numbers $\gamma_j \rightarrow
\infty$, $|j| > 1$, with $\gamma_{-j} = \gamma_j$, and for each $j$ a solution
$\Phi_\e^{j,\pm}(s)$ of $ \Lej \phi =0$ such that 
\[
 e^{\pm \gamma_j s}\Phi_\e^{j,\pm}(s),
\]
are {\it periodic} functions of $s$. In
particular, 
\begin{equation}
|\Phi_\e^{j,+}(s)| \le c e^{-\gamma_j s}, \qquad |\Phi_\e^{j,-}(s)| 
\le c e^{\gamma_j s}, \qquad \mbox{for all} \quad s \in {\RR}.
\label{eq:4.65}
\end{equation}
In fact, $\Phi_\e^{j,-}(s) = \Phi_{\e}^{j,+}(-s)$. Because of Corollary~\ref{cor:4.1}, 
it is natural to define $\gamma_0 = \gamma_{\pm 1} = 0$. 

While there are analogous conclusions were we to be using the independent 
variable $t$ instead of $s$, the values of these `indicial exponents' $\gamma_j = 
\gamma_j(\e)$ would behave in a less desirable way as $\e$ tends to
zero. 

\begin{proposition}
For any $\eta > 0$, there exists an $\e_0 > 0$ such that when $0 < \e \leq 
\e_0$, the 
numbers $\gamma_j$ satisfy $\gamma_j \geq \sqrt{4-\eta}$ for $|j| > 1$. 
\label{pr:4.15}
\end{proposition}
{\bf Proof:} The $\Phi_\e^{j,\pm}(s)$ are homogeneous solutions for
the ordinary differential operator $-\del_s^2 + Q_j$, where
$Q_j(s) = -\tau^2 \cosh 2\sigma(s) + j^2$. We are trying to estimate
the exponential growth rate of solutions for this operator. It follows from 
(\ref{eq:3.15}) that
\[
\tau^2 \leq \tau^2\cosh 2\sigma \leq 2-\tau^2.
\]
Thus, we see that $Q_j \geq j^2 -  2+\tau^2$. In particular, for $|j| \geq 3$,
$Q_j > 4$. In this case, the result is clear: $e^{\pm 2 s}$ are
supersolutions for the operator, and homogeneous solutions 
bounded by these supersolutions may be constructed by the shooting
method, as described next for the slightly more complicated case
when $j = \pm 2$. 

Let $w$ be the unique decreasing solution of $(-\del_s^2 + Q_2)w = 0$
with $w(0) = 1$. This solution may be constructed as a limit, as $s_1 
\rightarrow \infty$, of solutions $w_{s_1}$ of this equation on $[0,s_1]$ with 
$w_{s_1}(0) = 1, \ w_{s_1}(s_1) = 0$. Since $e^{-\sqrt{2+\tau^2}s}$ 
is a supersolution for this operator and dominates $w_{s_1}$ at the endpoints,
$s=0$ and $s=s_1$, it gives an upper bound for $w_{s_1}$ on the whole interval 
$0 \leq s \leq s_1$. Thus the limit as $s_1 \rightarrow \infty$, which
we call $w$, exists and is also bounded by this same function. Next, let 
$\bar{Q}_2 =  - (2/\cosh^2 s) + 4$, and define $\bar{w}$ to be the unique 
solution of $(-\del_s^2 + \bar{Q}_2)\bar{w} = 0$ with $\bar{w}(0) = 1$
which is decreasing on $[0,\infty)$. Since $\bar{Q}_2$ converges exponentially
 to $4$ as $|s| \rightarrow \infty$, we know that $\bar{w}(s) \sim C_+
e^{- 2s}$ as $s \rightarrow + \infty$; because it decreases as $s \rightarrow 
+\infty$ and because  this equation has no global bounded solutions, we deduce 
that $\bar{w}(s) \sim C_{-}e^{-2s}$ as $s \rightarrow -\infty$.

Next, for any fixed $s_0 > 0$, $Q_2$ converges in ${\cal C}^\infty$ to 
$\bar{Q}_2$ on the interval $[-s_0,s_0]$ as $\tau \rightarrow 0$; we claim
that $w \rightarrow \bar{w}$ uniformly, along with all its derivatives, on 
this interval too. To establish this, it suffices to show that the
Cauchy data of $\bar{w}$ at $s=0$ converges to that of $w$.
Indeed, using the equations satisfied by both $w$ and $\bar{w}$ we get 
\[
w \, \del_s^2 \bar{w} -\bar{w}\, \del_s^2 w = (\bar{Q}_2- Q_2)\,  w \, \bar{w}.
\]
An integration by parts leads to
\[
w_s (0)\, \bar{w}(0) - w(0)\, \bar{w}_s (0) =\int_0^{+\infty}(\bar{Q}_2- Q_2)
\, w \, \bar{w} ds.
\]
Thanks to Corollary \ref{cor:3.1} and also thanks to the fact that we already 
know that $\bar{w}(s) \sim C_+ e^{- 2s}$ as $s \rightarrow + \infty$ and 
that $|w| \leq c  e^{-\sqrt{2+\tau^2}s}$, we see that  the right hand side of 
this equality tends to $0$ as $\tau$ tends to $0$, and so
\[
\lim_{\tau \rightarrow 0} w_s (0) =\bar{w}_s (0).
\]
The claim then follows at once using Corollary \ref{cor:3.1} once again.
 
By the Bloch wave theoretic  construction of solutions of operators with 
coefficients periodic of period $S_\e /2$, we may write $w(s) = e^{-\gamma_2 s}
p(s)$, where $p(s+ S_\e /2) = e^{i\lambda}p(s)$ for some $\lambda \in {\RR}$.
 Since $p(s)$ is real valued and  strictly positive, 
actually $\lambda = 0$.  We wish to show that this exponent $\gamma_2$
converges to the exponential rate of decrease $-2$ corresponding to $\tau = 0$,
or more quantitatively, that for any fixed $\eta > 0$, we have  $\gamma_2 > 
\sqrt{4-\eta} \equiv \gamma$, so long as $\tau$ is sufficiently close to $0$. 

To estimate this exponential rate of decrease, it suffices to show
that for some fixed $s_0 > 0$, 
\begin{equation}
w( \ha S_\e -s_0) \le C \displaystyle e^{-\gamma \frac{S_\e}{2}} w(-s_0),
\label{eq:vit}
\end{equation}
where the constant $C$ is independent of $\tau$. We do this in two steps. 

First we  choose $s_0$ sufficiently large so that $2/\cosh^2 s_0 < \eta$;
also, using Proposition~\ref{pr:3.6}, we can assume that $ 2 \tau^2 \cosh^2 
\sigma \leq \eta$ on $[s_0, S_\e/2 -s_0]$ for all $\tau$ small enough.
We have noted that $\bar{w} \sim C_{\pm}e^{-2s}$ as $s \rightarrow \pm 
\infty$, and this implies in particular that, for $s >0$
\[
\bar{w}(s) \le   C e^{-4s}\bar{w}(-s),
\]
for some $C>0$ independent of $s$.
Moreover, increasing the value of $s_0$, if necessary, we may assume that 
\[
\bar{w}_s (s_0) \leq -(2-\eta/4) \bar{w} (s_0).
\]
By the uniform convergence of $w$ to $\bar{w}$, we deduce that, for $\tau$ 
small enough
\begin{equation}
w (s_0) \le  2  C e^{-4s_0} w (-s_0).
\label{eq:vvit}
\end{equation}
and also
\[
w_s (s_0) \le -(2-\eta /2)w(s_0).
\]
Next, on the interval $[s_0, S_\e /2 - s_0]$, the potential $Q_2$ is bounded 
from below  by $\gamma^2 \equiv 4- \eta$, and so we may use $\tilde{w}(s) = 
e^{-\gamma (s-s_0)}w(s_0)$ as a supersolution here. Indeed $w(s_0) = 
\tilde{w}(s_0)$ and 
\[
w_s (s_0) \le -(2-\eta /2)w(s_0) = -(2-\eta /2)\tilde{w}(s_0) 
=\frac{\sqrt{4-\eta}}{2} \tilde{w}_s (s_0) \leq \tilde{w}_s (s_0) ,
\]
provided $\tau$ is chosen small enough. We conclude that $w \le \tilde{w}$ 
on this whole interval. Thus, we get in particular the estimate
\[
w(\ha S_\e - s_0) \le \displaystyle e^{-\gamma (\ha S_\e - 2s_0)} w(s_0).
\]
Putting (\ref{eq:vvit}) together with this last inequality, we get that
\[
\displaystyle w(\ha S_\e - s_0) \le C e^{2s_0(\gamma -2)} e^{-\gamma \frac{S_\e}{2}}
w(-s_0) \le C e^{-\gamma\frac{S_\e}{2}}w(-s_0)).
\]
This gives the desired estimate, and the proof is complete.
\hfill $\Box$

\subsection{Mapping properties of the Jacobi operator} To fully analyze
the problem (\ref{eq:4.6}), we must study the mapping properties of the Jacobi 
operator $\Le$, both for fixed $\e > 0$ and uniformly down to $\e=0$. To 
state this result it is first necessary to define appropriate function spaces 
on which the Jacobi operator will act; these are exponentially weighted 
H\"older spaces. This is one of the main places where the difference between 
the independent variables $s$ and $t$ is seen: it is possible to obtain good 
mapping properties on spaces of this type, defined either in terms of the $s$ 
or $t$ variables, for fixed $\e$, but it is impossible to obtain
the uniform behaviour down to $\e=0$ when using $t$. Fortunately, this 
uniformity does occur when using $s$, and so from now on, unless saying
explicitly otherwise, this choice of independent variable will be used in the 
sequel. 

The definition of the weighted H\"older spaces is the natural one:
\begin{definition} Parametrize ${\RR}^+ \times S^1$ by the
variables $(s,\theta)$. For each 
$r \in {\mathbb N}$ and  $0<\alpha<1$ and $s \in {\R}$, let 
\[ 
|w|_{r,\al,[s, s +1]},
\]
denote the usual ${\cal C}^{r,\al}$ H\"older norm on the set $[s,s +1] \times
S^1$. Then for any $\mu \in {\RR}$ and $s_0 \in {\R}$, 
\[
\begin{array}{rlll}
{\cal C}^{r,\al}_\mu ([s_0, +\infty) \times S^1) & = & \left\{ w \in 
{\cal C}^{r,\al}_{\mathrm loc} ([s_0, +\infty) \times S^1) \quad 
\mbox{\rm and }\right.\\[2mm]
\, & \, & \left. \quad ||w||_{r,\al,\mu }= \sup_{s \geq s_0}e^{- \mu s}
|w|_{r,\al,[s,s + 1]} < \infty \right\}. 
\end{array}
\] 
\label{de:4.1} 
\end{definition} 
In particular, the function $e^{\mu s}$ is in ${\cal C}^{r,\alpha}_\mu
([s_0, +\infty) \times S^1)$. 

Recall now the splitting of $\Le$ into $\Le'$ and $\Le''$, corresponding
to the operator induced on the eigenspaces with $|j|\leq 1$ and
$|j| > 1$, respectively. $\Pi'$ and $\Pi''$ are the projectors 
onto the corresponding subspaces. This
will often be abbreviated by letting $\Pi' w = w'$, and so on. 
The main result of this section is the
\begin{proposition} 
Fix $\mu$ with $1 < \mu < 2$. Then there exists an $\e_0 >0$, depending only 
on $\mu$, such that whenever $\e \in (0, \e_0)$, there exists a unique 
solution $w \in {\mathcal C}^{2,\al}_{-\mu}([S_\e /8,\infty) \times S^1)$ 
of the problem 
\begin{equation} 
\left\{ 
\begin{array}{rlll} 
  \Le w & = & f \quad  & \mbox{in}\quad (S_\e /8,\infty) \times S^1\\[2mm] 
\Pi'' w & =  &  \phi'' \quad & \mbox{on}\quad \{S_\e /8\} \times S^1,
\end{array} 
\right. 
\label{eq:4.7} 
\end{equation} 
for $f \in {\cal C}^{0,\al}_{-\mu}([S_\e /8,\infty)\times S^1)$ and $\phi'' 
\in \Pi''\left( {\cal C}^{2,\al}(S^1)\right)$. The solution of 
the homogeneous Dirichlet 
problem, when $\phi'' = 0$, will be denoted $w = G_{\e}(f)$, while 
the Poisson operator, which gives the solution when $f=0$, will be 
denoted by $w=P_{\e}(\phi'')$. The linear maps
\[
\begin{array}{rlll}
G_{\e}:{\cal C}^{0,\al}_{-\mu} ([S_\e /8,\infty)\times S^1) & \longrightarrow 
& {\cal C}^{2,\al}_{-\mu}([S_\e /8,\infty)\times S^1)\\[2mm]
\e^{\mu /4} P_{\e}: \Pi''\left( {\cal C}^{2,\al}(S^1)\right)
 & \longrightarrow & {\cal C}^{2,\al}_{-\mu}([S_\e /8,\infty) \times S^1),
\end{array}
\]
are bounded uniformly for all $\e \in (0, \e_0)$.
\label{pr:4.3} 
\end{proposition} 
{\bf Proof:} The proof of the existence of $G_{\e}$ and $P_{\e}$
 and of their uniformity is accomplished in a number of steps. 
Solutions are constructed on each eigenspace of the Laplacian on $S^1$, and 
the cases where $|j|\leq 1$ must be  be treated somewhat differently than 
the others. 

We shall give the proof of this result in a slightly more general 
context where the boundary point $S_\e /8$ is replaced by $s_0$ arbitrarily 
chosen in ${\RR}$.

Fix $f$ and $\phi$ in the appropriate function spaces. We decompose  $w= w' 
+ w''$, $f= f'+ f''$, then we must solve 
\begin{equation} 
\left\{  
\begin{array}{rlll} 
L_{\e}' w' &   = & f'\qquad  &\mbox{for} \quad s >s_0 \\[2mm] 
L_{\e}'' w'' &  = & f'' \qquad & \mbox{for} \quad s >s_0 \\[2mm] 
w'' & = & \phi''  \qquad & \mbox{for} \quad s =s_0 .
\end{array} 
\right. 
\label{eq:4.8} 
\end{equation} 
Notice that, no boundary conditions are imposed on $w'$ at $s=s_0$.
We will also need to decompose
\[
w'(s, \theta) = \frac{1}{\sqrt \pi} w_{-1}(s) \sin \theta + \frac{1}{\sqrt{2\pi}} w_0(s) +
\frac{1}{\sqrt \pi} w_1 (s) \cos \theta,
\]
and 
\[
f'(s, \theta) = \frac{1}{\sqrt \pi} f_{-1}(s) \sin \theta
 + \frac{1}{\sqrt{2\pi}} f_0(s) + \frac{1}{\sqrt \pi} f_1 (s) \cos \theta.
\]

\medskip

{\bf Step 1:} We first consider the problem where  $\phi'' = 0$. Thus $f 
\in {\cal C}^{0,\al}_{-\mu}$, and multiplying by a suitable factor, 
we may assume that 
\[ 
||f'||_{0,\al,-\mu}+ ||f''||_{0,\al,-\mu} =1. 
\] 
In this step, we 
only consider the restriction of the problem to the high eigencomponents. 
We first show that for every $s_1 > s_0$ there is a unique solution of 
\begin{equation}
\left\{ 
\begin{array}{rllll} 
L_\e'' w''_{*}& = & f'' \quad & \mbox{in} \quad (s_0,s_1) \times S^{1}\\[2mm]  
w''_{*} & = & 0 \quad & \mbox{on} \quad \{s_0\} \times S^{1}\\[2mm]  
 w''_{*} & = & 0\quad & \mbox{on} \quad \{s_1\} \times S^{1}. 
\end{array} 
\right. 
\label{eq:ws}
\end{equation}
The existence of $w''_{*}$ follows from a standard variational argument 
using the energy functional 
\[ 
{\cal E}(w)=\int_{s_0}^{s_1}\int_{S^{1}}(|\del_s w|^2 + |\del_\theta w
|^2 - \tau^2 \cosh (2\sigma) \, |w|^2 + f'' w)\,ds\, d\theta. 
\] 
Using the fact that $\forall |j| >1$, we have  $j^2 - \tau^2 \cosh (2\sigma) 
>2$, we see that when we restrict the domain of ${\cal E}$ to the span of
the eigenfunctions  $\chi_j(\theta)$ with $|j|>1$, this functional is convex 
and proper, and the existence of a unique minimizer for it, which we denote by
$w''_{*}$, is then immediate. 
 
We claim that there exists a constant $\e_0>0$ and a constant 
$C = C(\mu)>0$, independent of $s_0 < s_1$ and $\e \in (0, \e_0)$, such that 
\[ 
\sup_{\theta \in S^{1}}\sup_{s\in [s_0,s_1]} e^{\mu s}|w''_{*}(s,\theta)| 
\leq C(\mu).
\] 
Assuming that the claim is already proven, we can choose a sequence $s_{1,i}$ 
tending to $+\infty$ and build $w''_{*,i}$ the corresponding solutions of 
(\ref{eq:ws}). The uniform bound above allows us to extract from the sequence 
$w''_{*,i}$ a subsequence which converges to a solution $w''$ of
\[
\left\{ 
\begin{array}{rllll} 
L_\e'' w''& = & f'' \quad & \mbox{in} \quad (s_0,+\infty) 
\times S^{1}\\[2mm]  
w'' & = & 0 \quad & \mbox{on} \quad \{s_0\} \times S^{1}  ,
\end{array} 
\right. 
\]
which satisfies
\[
\sup_{\theta \in S^{1}}\sup_{s\in [s_0,s_1]} e^{\mu s}|w''(s,\theta)| 
\leq C(\mu).
\] 
and then, by classical elliptic estimates, that
\[
||w''||_{2, \alpha, -\mu} \leq c(\mu) ,
\]
for some constant $c(\mu)>0$ independent of $s_0 \in {\R}$ and
 $\e \in (0, \e_0)$.

The claim is proved by contradiction. By assumption, we have
$e^{\mu s}|f''(s,\theta)|\leq 1$ for $s_0 \leq s \leq s_1$, $\theta \in S^1$. 
If the assertion were not true, then there would exist sequences of numbers
$s_{0,i}$, $s_{1,i}$, functions $f''_i$, Delaunay parameters $\e_i$
 and corresponding solutions ${w}''_{*,i}$ such that  
\[ 
A_i \equiv \sup_{\theta \in S^{1}}\sup_{s_{0,i} \le s \le s_{1,i}}
e^{\mu s}|w''_{*,i}(s,\theta)| \longrightarrow \infty,   
\] 
and 
\[
\sup_{\theta \in S^{1}}\sup_{s_{0,i} \le s \le s_{1,i}}
e^{\mu s}|f''_{i}(s,\theta)|\leq 1.  
\] 
Suppose that this maximum, for each $i$, is attained at some point
$(s_i,\theta_i)$, and define  
\[ 
\tw_i''(s,\theta) = A_i^{-1}e^{\mu s_i}w''_{s_{1,i}}(s+s_i, \theta), \qquad 
\tilde{f}''_i(s,\theta) = A_i^{-1}e^{\mu s_i}f''(s+s_i,\theta). 
\] 
Then
\[ 
\sup_{\theta \in S^{1}}\sup_{s_{0,i}-s_i \le s \le s_{1,i}-s_i}e^{\mu s}| 
{\tw}''_i(s,\theta)|= 1, 
\] 
and this supremum is attained on $\{0\}\times S^1$, while 
$\tilde{f}''_i \rightarrow 0$ in norm. Furthermore, 
\[ 
\Le'' \tw''_i   = \del_s^2 \tw''_i + \del_\theta^2 \tw''_i +
\tau_i^2 \cosh (2\sigma_i(\cdot +s_i) ) \tw''_i = \tilde{f}''_i 
\] 
on $[s_{0,i}-s_i, s_{1,i}-s_i]\times S^{1}$.  

Passing to a subsequence if necessary, we assume that $s_{0,i}-s_i$ converges 
to  $v_1 \in {\RR} \cup \{-\infty \}$ and $s_{1,i}-s_i$ converges to $v_2 
\in {\RR} \cup \{+\infty\}$. By using the result of Proposition~\ref{pr:3.7},
 the bounds above, as well as those 
provided by elliptic estimates, we can also assume that $\tw''_i$ converges, 
along with all its derivatives, over any compact subset of $(v_1,v_2)\times 
S^1$ (including endpoints if either is finite) to a function $\tw''$, which 
satisfies $e^{\mu s}|\tw''| \le 1$ over this set, is nonvanishing (because of 
the normalization of $\tw''_i$ at $s=0$), and which solves one of the 
following equations: 
\begin{equation}
\del_s^2 \tw'' + \del_\theta^2 \tw'' +\tau^2 \cosh (2\sigma (s +\bar{s}))
\tw'' =  0,
\label{eq:4.9}
\end{equation}
for some $\e \in (0, \e_0)$ and $\bar{s} \in {\RR}$,
\begin{equation}
\del_s^2 \tw'' +  \del_\theta^2 \tw''  = 0,
\label{eq:4.10}
\end{equation}
or  
\begin{equation}
\del_s^2 \tw'' + \del_\theta^2 \tw'' +\frac{2}{\cosh^2 (s+\bar{s})}\tw''= 0,
\qquad \mbox{ for some} \quad\bar{s}\in {\RR} ,
\label{eq:4.11}
\end{equation}
on $[v_1,v_2] \times S^{1}$. In addition, if either $v_1$ or $v_2$ is finite,
then $\tw''$ vanishes at that endpoint.

We must analyze a few cases, depending on the values of $v_1$ and  $v_2$ 
and which of the equations above is satisfied by $\tw''$. The goal
in each case is to show that $\tw''$ must, in fact, vanish
identically, which would be a contradiction. 

The point, in all cases, is that we wish to multiply the appropriate
equation for $\tw''$ by $\tw''$ and integrate by parts, to obtain
\[
\int_{v_1}^{v_2} \int_{S^1} |\del_s \tw''|^2 + |\del_\theta \tw''|^2
- \tau^2 \cosh (2\sigma (s + \bar{s}))\,  |\tw''|^2\,ds\,d\theta = 0, 
\]
\[
\int_{v_1}^{v_2}\int_{S^1} |\del_s \tw''|^2 + |\del_\theta \tw''|^2 \,ds\,
d\theta =0  
\]
or 
\[
\int_{v_1}^{v_2} \int_{S^1} |\del_s \tw''|^2 + |\del_\theta \tw''|^2
- \frac{2}{\cosh^2(s+\bar{s})} |\tw''|^2\,ds\,d\theta = 0.
\]
In each of these three cases we see that the integrand is positive, because 
we always have the inequality
\[
\int_{v_1}^{v_2} |\del_\theta \tw''|^2 \, ds \geq 4 \int_{v_1}^{v_2}
 |\tw''|^2 \, ds,
\]
and so we would conclude that $\tw'' \equiv 0$, which is a contradiction.

To make this argument work, it suffices to show that the boundary
terms in the integration by parts vanish. When either $v_1$ or $v_2$
is finite, this is immediate from the Dirichlet conditions
at that boundary, so it remains to show that if either $v_1$ or $v_2$ 
is infinite, then $\tw''$ decays exponentially in that direction.
Any unbounded solution of (\ref{eq:4.10}) on a half-line must grow 
at least at the rate $e^{2|s|}$, which would violate the condition
$e^{\mu s} |\tw''| \le 1$, so we see that $\tw''$ must decrease exponentially 
in this case. The same argument works when $\tw''$ satisfies (\ref{eq:4.11}) 
because solutions of that equation have the same asymptotic rates
of growth or decay as solutions of (\ref{eq:4.10}). Finally, if $\tw''$ 
satisfies  (\ref{eq:4.9}), we first choose $\eta >0$ such that  
$\sqrt{4-\eta} > \mu$, then, we apply Proposition~\ref{pr:4.15}, which states
that any unbounded solution must grow at least at the rate
$e^{\sqrt{4-\eta} |s|}$ provided $\e$ is less than, say, $\e_0$. Therefore,
we can eliminate the possibility of exponential growth. This ends the proof 
of the claim.

\medskip

{\bf Step 2:} We now consider the cases when $|j| \leq 1$. 
The argument when $j = \pm 1$ is almost identical to the one for $j=0$, so
we shall just consider the latter case, commenting on the
end on the very minor changes that need to be made. 
Thus, recalling that we are no longer requiring any boundary
conditions, we wish to find a solution to the problem 
\begin{equation} 
L_{\e,0} w_0 \equiv \del_s^2 w_0 +\tau^2 \cosh (2 \sigma) w_0 = f_0 \qquad 
\mbox{in} \quad [s_0, \infty),
\label{eq:4.12} 
\end{equation} 
with the desired decay property at infinity.
We find this solution again as a limit of functions $w_{*}$ solutions of 
$L_{\e, 0} w_* =f_0$ on $[s_0, s_1)$, where now $w_{*}(s_1) = \del_{s}
w_{*}(s_1) = 0$. For convenience,
we choose a ${\cal C}^{0,\al}$ extension of $f_0$, vanishing when
$s < s_0 - 1$, say, and consider the solution $w_{*}$ for this
extended right hand side, now defined on $(-\infty, s_1]$. 

As in Step 1, we claim that there exists a constant $C = C(\mu)$, independent 
of $s_0, s_1$ and $\e$, such that 
\[ 
\sup_{s\in (-\infty,s_1]}e^{\mu s}|w_{*}(s)|\leq C. 
\] 
Once this claim is proved, the arguments of the proof are identical to those
in Step 1, so we shall omit them.

Again this is proved by contradiction. First, note that when $s < s_0 - 1$,
$w_{s_1}$ is a linear combination of the Jacobi fields $\Phi_{\e}^{0,\pm}$,
hence is at most linearly growing. If the assertion were false, there would 
exist sequences $f_{0,i}$, $s_{0,i}$, $s_{1,i}$, $\e_i$, and $w_{*,i}$ 
such that 
\[
A_i \equiv \sup_{s\in (-\infty , s_{1,i}]}e^{\mu s}|w_{*,i}| 
\longrightarrow \infty,
\]
and 
\[
\sup_{s\in (-\infty , s_{1,i}]}e^{\mu s}|f_{0,i}| \leq 1.
\]
If this maximum is attained at  $(s_i,\theta_i)$, $s_i \in (-\infty , 
s_{1,i})$, then we rescale the functions and translate the independent 
variable by $s_i$ to obtain a solution of 
\[ 
\frac{d^2{\tw}_i}{ds^2} +\tau_{i}^2 \cosh (2\sigma_i){\tw}_i = 
\tilde{f}_{0,i}, 
\] 
in $(-\infty ,s_{1,i} - s_i]$ which satisfies \[
\sup_{s\in (-\infty,s_{1,i}-s_i]}e^{\mu s} |{\tw}_i (s)| =1,
\] 
while $\tilde{f}_{0,i}$ tends to zero in norm. 
 
Passing to a subsequence, we obtain in the limit a nontrivial solution 
${\tw}$ of one the following equations:
\[ 
\frac{d^2{\tw}}{ds^2} +\tau^2 \cosh (2\sigma (s +\bar{s})) \tw = 0, \quad 
\mbox{for some} \quad \e\in (0, \e_0) \quad \mbox{and} \quad\bar{s} \in {\RR}, 
\] 
\[ 
\frac{d^2{\tw}}{ds^2} +\frac{2}{\cosh^2 (s+\bar{s})} \tw = 0,\quad \mbox{for 
some} \quad \bar{s}\in {\RR},
\] 
or 
\[ 
\frac{d^2{\tw}}{ds^2}= 0, 
\] 
over some interval $(-\infty ,v]$, and in each case, $|\tw (s)|\leq 
e^{- \mu s}$ in $(-\infty, v]$. 

Clearly $v$ cannot be finite, because if it were then  $\tw$ would have
to satisfy $\tw (v)=\del_s \tilde{w} (v)=0$, which would imply that it
would vanish identically. 

Now, for each of the three equations we know that there are no exponentially
decreasing solutions; for the second and third equations this is obvious,
while for the first it follows because we know the family of solutions
explicitly.  However, since we know that $\tw$ does decay exponentially
as $s \rightarrow \infty$, we again would have to conclude that it
vanishes identically, and this is a contradiction. 

When $j = \pm 1$, the changes that need to be made in this argument 
are minor. For example, when $j=1$, for all $s_1 >s_0$, the solution $w_{*}$ 
is defined as before to be the solution of 
\[
L_{\e,1} w_* \equiv \del_s^2 w_* - w_* + \tau^2 \cosh (2 \sigma) w_* = f_{1} 
\qquad  \mbox{in} \quad [s_0, s_1),
\]
which satisfies $w_* (s_1)=\del_s w_*(s_1)=0$ and where $f_1$ has been 
extended by $0$ in  $(-\infty, s_0-1]$.

And to establish its uniform bound, we proceed by contradiction. In this
case, however, the limiting equations are now
\[ 
\frac{d^2{\tw}}{ds^2} -\tw +\tau^2 \cosh (2\sigma (s + \bar{s})) \, \tw = 0,
 \quad \mbox{for some} \quad  \e \in (0, \e_0) \quad \mbox{and} \quad \bar{s} 
\in {\RR},
\] 
\[ 
\frac{d^2{\tw}}{ds^2} -\tw +\frac{2}{\cosh^2 (s+\bar{s})}\tw = 0,\quad 
\mbox{for some} \quad  \bar{s}\in {\RR},
\] 
or 
\[ 
\frac{d^2{\tw}}{ds^2} -\tw =0, 
\] 
on $(-\infty, v] \times S^1$, with boundary condition $\tw(v)=\del_s \tw
(v)=0$ if $v$ is finite, and where $|\tw (s)|\leq e^{- \mu s}$ for all 
$s \in (-\infty, v]$. 

Once again, $v$ cannot be finite, but now the equations do admit
exponentially decreasing solutions at $\pm \infty$. However, all such solutions
decay no faster than $e^{-s}$, whereas we have assumed that $1 < \mu < 2$,
so once again we obtain a contradiction. 

\medskip

{\bf Step 3:}  Finally consider the problem when $f=0$ and $\phi'' \neq 0$. 
We may as well assume that $\|\phi''\|_{2,\al} = 1$. Let $\eta(s)$ be a smooth 
cutoff function equal to $1$ for $s\leq 0$ and vanishing for $s \geq 1$. Then 
\[
\Le'' w = 0, \qquad w (s_0,\theta) = \phi'' (\theta), 
\]
is equivalent to 
\[
\Le'' \bar{w} = - \Le'' (\eta (s-s_0) \phi''(\theta)) , \qquad w 
(s_0,\theta) = 0, 
\]
which has already been solved in Step 1. Moreover, since 
\[
||\eta \phi''||_{0, \alpha, -\mu} \leq c \, e^{\mu s_0},
\]
it follows from Step 1 that
\[
||w||_{2, \alpha, -\mu} \leq c \, e^{\mu s_0},
\]
as we wished. This completes the proof in all cases.  \hfill $\Box$

\begin{corollary}
Fix $1< \mu < 2$. Then there exists a constant $c>0$ and 
an $\e_0 >0$, depending only on $\mu$, such that for $0<\e<\e_0$, we have
\[
||(P_{\e} - P_{0})(\phi'')||_{2,\al,-\mu} \leq c \, \e^{- \mu /4} 
\left(\e^{1/2} + \e^{(6-3\mu)/4}\right)||\phi''||_{2,\al}.
\]
Here, if $\phi'' \in \Pi''\left({\cal C}^{2,\al}(S^1)\right)$, the function
$P_{0}(\phi'')$ is the unique solution in ${\cal C}^{2,\al}_{-2}
([S_\e /8,\infty)\times S^1)$ of the problem
\begin{equation} 
\left\{ 
\begin{array}{rlll} 
\Delta w & = & 0 \quad & \mbox{in}\quad [S_\e /8,\infty) \times S^1\\[2mm] 
 w  & =  &   \phi''\quad & \mbox{on}\quad \{S_\e /8\} \times S^1.
\end{array} 
\right. 
\label{eq:lap}
\end{equation} 
\label{cor:11}
\end{corollary}
{\bf Proof:} 
Write $w_\e = P_{\e}\phi''$ and $w_0 = P_{0}\phi''$. 
If $w_\e = w_0 + h$, then $L_\e h = -(\tau^2 \cosh (2\sigma)  )\, w_0$
and $\Pi'' h(S_\e /8,\theta) = 0$, and so $h = -G_{\e}( \tau^2
\cosh (2\sigma)\,  w_0)$. We first estimate 
\[
||h||_{2,\al,-\mu} \leq c\, || \tau^2 \cosh (2\sigma) \,w_0||_{0,\al,-\mu} .
\]
Using
\[
||w_0||_{0,\al,[s,s+1]} \leq c\,e^{-2(s-S_\e /8)}||\phi''||_{0,\al} \leq
c\,\e^{-\ha} e^{-2s}||\phi''||_{0,\al},
\] 
we bound this by 
\[
c\,\e^{- \ha}\left(\sup_{s \geq S_\e /8} e^{(\mu-2)s}
||\tau^2 \cosh 2\sigma ||_{0,\al,[s,s+1]}\right)\,||\phi''||_{0,\al}.
\]
When $S_\e/8 \leq s \leq  3 S_\e /8$, we know from Proposition~\ref{pr:af} 
that
\[
\tau^2\cosh (2 \sigma) + |\del_s\tau^2\cosh ( 2 \sigma)|  \leq \e^{\ha}.
\]
Therefore
\[
\e^{-\ha}e^{(\mu-2)s}||\tau^2 \cosh 2\sigma ||_{0,\al,[s,s+1]} \leq c\,
\e^{(\mu-2) S_\e /8} = c\,\e^\ha \, \e^{ - \mu /4}.
\]
Next, when $ 3 S_\e /8 \leq s \leq S_\e /2$, we know, still from 
Proposition~\ref{pr:af}, that we may estimate 
\[
\tau^2 \cosh 2\sigma + |\del_s\tau^2 \cosh 2\sigma | \leq c\,\e^2 e^{2s},
\]
hence
\[
\e^{-\ha}e^{(\mu-2)s}||\tau^2 \cosh 2\sigma||_{0,\al,[s,s+1]}
\leq c\,\e^{-\ha + 2 - \mu} = c\,\e^{(6-3\mu)/4}\, e^{- \mu/4}.
\]
Finally, for $S_\e /2 \leq s$ we use the fact that $\tau^2 \cosh 2\sigma 
\leq 2$, and proceed as before. This proves the Corollary. \hfill $\Box$
 
\section{CMC surfaces near to a half Delaunay surface}

In this section we construct by perturbation methods the
full space of CMC surfaces near to a fixed (half) Delaunay surface
$\D$ of necksize $\e$, as usual controlling the behaviour
as $\e \rightarrow 0$. Assume that $\D$ has the parametrization
\[
\x(s,\theta) = (\tau e^{\sigma (s)} \cos \theta,\tau e^{\sigma (s)} 
\sin \theta, k(s)),
\]
where $\tau \in (0,1)$, $\sigma$ and $k$ are as in 
Proposition~\ref{pr:3.1}. The unit normal at $\x(s,\theta)$ is defined to be
\[
\nu (s, \theta) = (-\tau \cosh \sigma(s) \cos \theta,
-\tau \cosh \sigma (s) \sin \theta, \sigma_s (s)).
\]
Therefore, surfaces which may be written as normal graphs over $\D$
admit the parametrization 
\[
\x_w(s, \theta) = \x (s, \theta) + w(s, \theta) \nu (s, \theta),
\]
for some sufficiently small function $w$ on $\D$. We denote by ${\D}_w$ the
surface obtained in this way. The components of its metric tensor are  
\[
E_w = \tau^2\left(e^\sigma - \sinh \sigma \, w\right)^2 + w_s^2,
\qquad F_w = w_s w_\theta, 
\]
and 
\[ 
G_w = \tau^2 \left( e^\sigma - \cosh \sigma \,w \right)^2 + w_\theta^2.
\]
The components of the second fundamental form are considerably 
less simple. In computing the following, we use that we always have 
the bounds  $\tau e^\sigma 
\leq 2$, $\tau^2 \sinh^2 \sigma \leq 1$ and $\sigma^2_s \leq 1$.
After substantial work, we find that 
\[
\sqrt{E_w G_w -F_w^2}\, L_w = \tau^3 e^{3\sigma} \left(
\tau \sinh \sigma + P_1 \left(\frac{w}{\tau e^{\sigma}}, 
\frac{\nabla w}{\tau e^{\sigma}}, \frac{\nabla^2 w}{\tau e^{\sigma}}
\right)\right),
\]
where $P_1$ is some polynomial (of degree at most $3$) without any constant
 term, the coefficients 
of which are functions of $s$ and such that they and their derivatives 
are bounded uniformly in $s$ and $\e$. In a similar manner we derive that
\[
\sqrt{E_w G_w -F_w^2}\, M_w = \tau^3 e^{3\sigma} 
P_2 \left(\frac{w}{\tau e^{\sigma}}, \frac{\nabla w}{\tau e^{\sigma}}, 
\frac{\nabla^2 w}{\tau e^{\sigma}}\right)
\]
and 
\[
\sqrt{E_w G_w -F_w^2} \, N_w =  \tau^3 e^{3\sigma} \left( \tau \cosh \sigma + 
P_3 \left(\frac{w}{\tau e^{\sigma}}, \frac{\nabla w}{\tau 
e^{\sigma}}, \frac{\nabla^2 w}{\tau e^{\sigma}}\right)\right),
\]
where $P_2$ and $P_3$ have the same properties as $P_1$. 

The equation that ${\D}_w$ has mean curvature $1$ is
\begin{equation}
L_w G_w - 2 M_w F_w + N_w E_w - ( E_w G_w -F_w^2) = 0.
\label{eq:5.1}
\end{equation}
This is a rather complicated nonlinear elliptic equation for $w$
which we shall not write out in full. Notice that it is satisfied
when $w=0$. Using the previous formula for the coefficients of the first 
and second fundamental forms, we find that its Taylor expansion about $w=0$ is
\begin{equation}
L_\e w = \tau e^{\sigma}  Q \left( \frac{w}{\tau e^{\sigma}}, 
\frac{\nabla w}{\tau e^{\sigma}}, \frac{ \nabla^2 w}{\tau e^{\sigma}} \right),
\label{eq:5.2}
\end{equation}
where 
\[
L_\e  w = w_{ss} + w_{\theta \theta}+ \tau^2 \cosh ( 2 \sigma) w,
\]
and $Q$ is again a polynomial (now of higher order) without any constant
 or linear terms, the coefficients of 
which have partial derivatives bounded uniformly in $s$ and $\e$.
We also write, for brevity,
\[
{\Q}(w)\equiv \tau e^{\sigma} Q \left( \frac{w}{\tau e^{\sigma}}, 
\frac{\nabla w}{\tau e^{\sigma}}, \frac{ \nabla^2 w}{\tau e^{\sigma}}\right).
\]

Given $\phi''\in \Pi'' \left({\cal C}^{2,\al}(S^1)\right)$, we would like 
to solve 
the boundary value problem
\begin{equation}
\left\{ \begin{array}{rlll}
L_\e w & = & {\cal Q} (w)   \qquad & \mbox{in} \qquad [S_\e /8, +\infty) 
\times S^1 \\[2mm]
\Pi'' w & = &  \phi'' \qquad & \mbox{on}  \qquad \{S_\e /8 \} \times S^1 .
\end{array}
\right.
\label{eq:5.4}
\end{equation}
Let $w_\e$ be the unique solution in 
${\mathcal C}^{2,\al}_{-\mu}([S_\e /8, +\infty)  \times S^1)$, $1<\mu<2$, of 
\[
\left\{ \begin{array}{rlll}
L_\e w_\e &  =  & 0  \qquad &\mbox{in} \qquad [S_\e /8, +\infty) 
\times S^1 \\[2mm]
\Pi'' w_\e  & = & \phi''\qquad &  \mbox{on} \qquad \{S_\e /8 \} \times S^1 ,
\end{array}
\right.
\]
which is given by Proposition~\ref{pr:4.3}. 
Setting $w = w_\e + v$, then we would like to find $v \in 
{\mathcal C}^{2,\al}_{-\mu}([S_\e /8, +\infty)  \times S^1)$ such that 
\[
\left\{ \begin{array}{rlll}
 L_\e v & = &  {\cal Q} (w_\e + v) & \qquad \mbox{in} \qquad 
[S_\e /8, +\infty) \times S^1 \\[2mm]
\Pi'' v  & = & 0 & \qquad \mbox{on} \qquad \{S_\e /8\} \times S^1.
\end{array}
\right.
\]
Notice that, it is sufficient to find a fixed point of the mapping
\begin{equation}
{\cal K}(v) = G_{\e}\Q(w_\e + v),
\label{eq:fpdel}
\end{equation}
at least when $\e$ is sufficiently small.
\begin{proposition} 
There exists a constant $c_0 >0$ such that if $||\phi''||_{2,\al} 
\leq c_0\,\e^{3/4}$, then 
\[
||G_{\e}({\Q}(w_\e ))||_{2,\al,-\mu}\leq c \,\e^{-3/4} 
||\phi''||^2_{2,\al} 
\]
and
\[
||G_{\e}({\Q}(w_\e+v_2)-{\Q}(w_\e +v_1))||_{2,\al,-\mu}
\leq \ha ||v_2-v_1||_{2,\al,-\mu},
\]
for all $v_1, v_2$ in $B_{c_0} \equiv \{v: ||v||_{2,\al,-\mu}\leq 
c_0\,\e^{(3-\mu)/4}\}$. Thus, ${\mathcal K}$ is a contraction mapping
on the ball $B_{c_0}$ into itself. Consequently, 
${\mathcal K}$ has a unique fixed point $v$ in this ball.
\label{pr:fixptdel}
\end{proposition}
{\bf Proof:} We shall use that 
\begin{equation}
||w_\e||_{2,\al,[s,s+1]} \leq c\, e^{\mu (S_\e /8-s)}||\phi''||_{2,\al}
\leq c\,c_0\,\e^{3/4}\,  e^{\mu(S_\e /8 -s)},
\label{eq:est1}
\end{equation}
and also that
\begin{equation}
||v||_{2,\al,[s,s+1]} \leq c_0\, \e^{3/4} e^{\mu(S_\e /8 - s)},
\label{eq:est2}
\end{equation}
for $v \in B_{c_0}$.

First consider $s$ in the range $[S_\e /8, 7 S_\e /8]$. Here, from 
Proposition~\ref{pr:af}, we get
\[
\tau e^{\sigma(s)}  \geq c \e^{3/4}.
\] 
Together with the fact that all derivatives of $\sigma$ are bounded, this gives
\[
\left\| (\tau e^\sigma) \left(\frac{w_\e}{\tau e^\sigma}
\right)^j\right\|_{2,\al,[s,s+1]} \leq c\, \e^{-3(j-1)/4}
||w_\e||_{2,\al,[s,s+1]}^j \leq c\, c_0^{j-2} e^{j\mu (S_\e /8 - s)}\, 
\e^{-3/4} \, ||\phi''||^2_{2,\al}
\]
for any integer $j \geq 2$. Thus, we already have obtained that
\[
e^{\mu s}||{\Q}(w_\e)||_{0,\al,[s,s+1]}
\leq c\,  e^{2\mu(S_\e /8-s) + \mu s}\, \e^{-3/4} ||\phi''||^2_{2,\al}
\leq c\, c_0^2\,\e^{(3-\mu)/4}. 
\]
Similarly
\[
e^{\mu s}||{\cal Q}(w_\e+v_1)-{\cal Q}(w_\e+v_2))||_{0,\al,[s,s+1]}
\]
can be estimated by the sum of products of 
$e^{\mu s}||v_1 - v_2||_{2,\al,[s,s+1]}$ with various terms of the form 
\[
||(w_\e + v_1)^j(w_\e + v_2)^{j'} (\tau e^{\sigma})^{-j-j'}||_{0,\al,[s,s+1]},
\qquad \mbox{\rm where}\ j+j' \geq 1.
\]
Each of these can be bounded by $c\,(c_0\,\e^{3/4})^{j+j'}\,
\e^{-(3/4)(i+j)}e^{(j + j')\mu (S_\e /8 - s)} \leq c\,c_0^{j+j'}$,
and so can be made as small as desired. 

For $s \geq 3 S_\e /4$, we have $s - S_\e /8 \geq 3 S_\e /4  \geq c - 
(3/2)\log \e$, and we will simply use the fact that
\[
\tau e^{\sigma(s)} \geq  \tau e^{\sigma(0)} = \e.
\]
Arguing as before, we get first that 
\[
e^{\mu s}||{\cal Q}(w_\e )||_{0,\al,[s,s+1]}
\leq c e^{\mu s} \e^{-1} e^{2\mu (S_\e /8 - s)} ||\phi''||^2_{2,\al}
\leq c \e^{-\mu /4 }(\e^{-1/4} e^{\mu (S_\e /8 - s)})\e^{-3/4}
 ||\phi''||^2_{2,\al}
\]
\[
\leq c_0\,\e^{(3-4\mu)/4}\, (c\,c_0\,\e^{-(1/4) + 3\mu/2}),
\]
and the quantity on the right in parentheses can be made as small as desired
when $\e$ is sufficiently small. Furthermore
\[
e^{\mu s}||{\cal Q}(w_\e+v_1)-{\cal Q}(w_\e+v_2))||_{0,\al,[s,s+1]}
\leq c\,c_0 \e^{-1} \e^{3/4}e^{\mu(S_\e /8-s)}
e^{\mu s}||v_2-v_1||_{2,\al,[s,s+1]}
\]
\[
\leq (c\,c_0\,\e^{-1/4 + 3\mu/2})||v_2-v_1||_{2,\al,-\mu},
\]
and again the coefficient can be made as small as desired when $\e$ 
is chosen small enough. 

Putting the estimates in these two domains together, and using that 
 $G_{\e}$ is bounded, we have now checked all the conditions necessary 
to ensure that ${\mathcal K}$ is a contraction mapping. Therefore there is
a unique element $v \in B_{c_0}$ such that ${\mathcal K}(v) = v$, and
the proof is complete. \hfill $\Box$

Examining this proof more carefully, we also obtain the 
\begin{corollary}
There exists a constant $c_0 >0$ and an $\e_0>0$ such that, for all $\e
\in (0, \e_0)$ and  for any $\phi'' \in \Pi''\left( 
{\cal C}^{2,\al}(S^1)\right)$ with $||\phi''||_{2,\al}\leq c_0\,\e^{3/4}$, 
the problem (\ref{eq:5.4}) has a unique solution $w$. The mapping 
\[
\Pi''\left( {\cal C}^{2,\al}(S^1)\right) \ni \phi'' \longrightarrow w \in 
{\cal C}^{2,\al}_{-\mu}([S_\e /8,\infty) \times S^1),
\]
is continuous and the solution $w$ satisfies the estimates
\begin{equation}
||w ||_{2,\al,-\mu}\leq c\, \e^{-\mu /4}(||\phi''||_{2,\al} + \e^{-3/4}
||\phi''||^2_{2,\al}),
\label{delrefest1}
\end{equation}
and
\begin{equation}
|| (w-\Pi''w ) (S_\e /8,\cdot)||_{2,\al}+||\del_s (w -
\Pi'' w) (S_\e /8,\cdot)||_{1,\al} \leq c \e^{-3/4}||\phi''||^2_{2,\al}.
\label{delrefest2}
\end{equation}
Finally, if $w_0 = P_{0}(\phi'') \in {\mathcal C}^{2,\al}_{-\mu} 
([S_\e /8, +\infty) \times S^1)$ as in Corollary~\ref{cor:11}, then 
\begin{equation}
||w -w_0 ||_{2,\al,-\mu} \leq c\, \e^{-\mu /4} \left( (\e^{1/2} +
\e^{(6-3\mu)/4})||\phi''||_{2,\al}+\e^{-3/4}||\phi''||^2_{2,\al}\right).
\label{delrefest3}
\end{equation}
\label{cor:1111}
\end{corollary}
{\bf Proof:} The solution $w$ is a sum $w_\e + v$ and we already know
that $||w_\e||_{2,\al,-\mu} \leq c \, ||\phi''||_{2,\al}$. For fixed $\phi''$, 
the map ${\mathcal K}$ is actually a contraction on the balls of radius a 
constant times $\e^{-(\mu+ 3)/4}||\phi''||^2_{2,\al}$, and so the norm of 
$v$ is at most this large. And this gives  (\ref{delrefest1}). The second 
estimate (\ref{delrefest2}) follows by evaluting at $s=S_\e /8$. Finally, for 
(\ref{delrefest3}), we write
\[
||w - w_0||_{2,\al,-\mu} \leq ||w_\e - w_0||_{2,\al,-\mu} + ||v||_{2,\al,-\mu},
\]
and use Corollary~\ref{cor:11}. \hfill $\Box$

\section{$k$-noids}

The second type of component in our construction of CMC surfaces are
a somewhat restricted class of minimal surfaces of finite total curvature with 
$k$ ends, or as we shall call them, $k$-noids. In this brief section we
discuss some of the global and asymptotic aspects of the geometry and topology
of $k$-noids, and in the next, discuss Jacobi operators on these surfaces 
and their compact truncations.

It is well-known that any $k$-noid $\Sigma$ has finite topology, and in fact 
is conformally equivalent to the complement of a finite number of points in a 
compact Riemann surface $\bar{\Sigma}$, i.e. $\Sigma = \bar{\Sigma}
\setminus \{p_1, \ldots, p_k\}$.  As in the introduction, we denote the space 
of $k$-noids of genus $g$ by ${\cal H}_{g,k}$. When $g > 0$, ${\cal H}_{g,1}$ 
and ${\cal H}_{g,2}$ are empty, while ${\cal H}_{0,1}$  and ${\cal H}_{0,2}$ 
contain only the plane and catenoid, respectively. 
The standard catenoid ${\cal C}_1$ is a surface of revolution, given in 
cylindrical coordinates by the parametrization
\[
\x(s,\theta) = (\cosh s \cos \theta, \cosh s \sin \theta, s).
\]
This is a conformal parametrization, and the unit normal, metric tensor and 
second fundamental forms are given by 
\begin{eqnarray*}
\nu(s,\theta) = \frac{1}{\cosh s} (- \cos \theta, -\sin \theta, \sinh s ),\\
g = \cosh^2 s \, ( ds^2 + d\theta^2), \qquad A = -ds^2 + d\theta^2.
\end{eqnarray*}
In particular, the mean curvature vanishes, and the catenoid is minimal.

We shall be discussing the space of moduli of $k$-noids. Just as with 
CMC surfaces, it is possible to determine the moduli space explicitly
in the simplest case, when $k=2$. In fact, the only complete minimal
surfaces in $\RR^3$ with two ends are images of the standard catenoid
${\cal C}_1$ by rigid motions and homotheties. While we shall 
frequently not distinguish between ${\cal C}_1$ and its translates
or rotations, it will be important to keep track of the homothety
factor. The dilation of ${\cal C}_1$ by the factor $a$ will be denoted 
${\cal C}_a$, and has the parametrization
\[
\x^{(a)}(s,\theta) = (a\cosh s \cos \theta, a\cosh s \sin \theta,a s).
\]
Thus any element of ${\cal H}_{0,2}$ is given as a rigid motion of some 
${\cal C}_a$. The metric tensor and second fundamental forms for this 
parametrization are
\begin{equation}
g_a = a^2 \cosh^2 s \,\big( ds^2 + d\theta^2\big), \qquad
A_a = a (- ds^2 +d\theta^2 ).
\label{eq:6.2}
\end{equation}

The plane and catenoid provide the asymptotic models for the ends 
of any $k$-noid: the basic structure theorem for $k$-noids states 
that an end of any $k$-noid may be written as a normal graph of a 
decaying function over an end of some suitably translated, rotated 
 plane or dilated catenoid. The corresponding ends will then be 
referred to as planar or catenoidal. Only $k$-noids with all ends 
catenoidal will be used in our construction; henceforth this will
always be assumed. 

Using this asymptotics theorem, we may assign a dilation, or
weight, parameter $a_\ell$ to each end $E_\ell$ of $\Sigma \in 
{\cal H}_{g,k}$, $\ell = 1, \ldots, k$, signifying that that end 
is the normal graph over (some translated and rotated copy of) 
${\cal C}_{a_\ell}$. This is analogous to the necksize parameters of 
the ends of CMC surfaces. This defines, at least in neighbourhoods
of the moduli space where some ordering of the ends is fixed, a map 
${\cal H}_{g,k} \rightarrow {\RR}^k$.

Fix $\Sigma \in {\cal H}_{g,k}$. We describe the parametrization of the ends
more carefully. Assume that $\Sigma$ has been rotated and translated
so that the the end $E_\ell$ is asymptotic to the model ${\cal C}_{a_\ell}$. 
By definition, there is a function $w$, defined on ${\cal C}_{a_\ell} \cap 
\{s \ge s_\ell\}$, such that $E_\ell$ is parametrized by 
\[
\x_w(s,\theta) \equiv \x_{a_j}(s,\theta) + w(s,\theta)\nu(s,\theta), \ 
s \geq s_\ell.
\]
This gives a canonical cylindrical coordinate system $(s,\theta)$ on 
$E_\ell$, which we will always use. The function $w$ is assumed {\it a 
priori} only to decay, but in fact admits an asymptotic expansion
\[
w(s,\theta) \sim \sum_{|j| > 1} a_j \chi_j(\theta) e^{-j s}, \qquad
\mbox{as} \quad s \rightarrow \infty.
\]

\section{The Jacobi operator on $k$-noids} 
Continuing our treatment of analysis on $k$-noids paralleling
that on Delaunay surfaces, we now consider the Jacobi operator 
$L = L_\Sigma$, which is the linearization of the mean curvature operator 
${\cal M}$ over $\Sigma$. This is
\[
L_\Sigma = \Delta_{\Sigma} + |A_\Sigma|^2,
\]
where the term of order zero is the squared norm of the second fundamental 
form of $\Sigma$. 

\subsection{Mapping properties of $L$ and Jacobi fields}
Just as for Delaunay surfaces, we require detailed knowledge of the mapping 
properties of $L$, first over all of $\Sigma$, then in a later section
for the Dirichlet problem on certain (deformations of) compact 
truncations $\Sigma_\e$ of $\Sigma$, and finally uniformly as 
$\e \rightarrow 0$. 

The analysis for $L$ over the complete surface $\Sigma$ is 
based on the fact that the ends have good asymptotic models. In 
fact, using the canonical cylindrical coordinates on each end $E_\ell$, we
see that the Jacobi operator for the model catenoid there is 
\begin{equation}
L_{a_\ell} = a_\ell^{-2} \cosh^{-2}s \left( \del_s^2 + \del_\theta^2
+ \frac{2}{\cosh^2 s}\right),
\label{eq:6.3}
\end{equation}
and so the true Jacobi operator is equal, as $s \rightarrow \infty$
in $E_\ell$, to the sum of this model operator and a correction
term, which is a second order operator each coefficient of which
decays at least like $e^{-4s}$. 
 
We let $L$ act on the weighted H\"older spaces ${\cal C}^{r,\al}_\mu 
(\Sigma)$, 
where $\phi$ is in this space if it is locally in ${\cal C}^{r,\al}(\Sigma)$ 
and on each end may be written as $e^{ \mu  s}\psi$ where $\psi \in 
{\cal C}^{r,\al}({\RR}^+_s \times S^1_\theta)$. 
 The basic mapping properties for $L$ are summarized in the
\begin{proposition} The operator
\[
L: {\cal C}^{2,\al}_{\mu}(\Sigma) \longrightarrow {\cal C}^{0,\al}_{\mu-2}
(\Sigma),
\]
is Fredholm provided $\mu \notin {\mathbb Z}$. In addition, $L$
is surjective on ${\cal C}^{2,\al}_\mu(\Sigma)$ if and only if
it is injective on ${\cal C}^{2,\al}_{-\mu}(\Sigma)$. 
\label{pr:6.1}
\end{proposition}

The drop of two in the weight parameter comes from the factor 
$(\cosh(a_{\ell}s))^{-2}$ in the expression for $L$ on $E_\ell$.
This sort of result is fairly standard by now; it may be proved by
constructing local parametrices, or solution operators, for the
model operators on each of the ends $E_\ell$ using explicit
ODE techniques on each of the cross-sectional eigenspaces, joining
these to a parametrix for the interior compact region, and finally
using standard perturbation techniques and Fredholm theory.

This Proposition leads naturally to the issue of determining the 
values of the weight parameters $\mu$ for which $L$ is surjective or 
injective. Although for a given $k$-noid $\Sigma$ this may be quite 
difficult to determine, the following condition is essential for the 
moduli space theory: 
\begin{definition}
A $k$-noid $\Sigma$ is called nondegenerate if its Jacobi
operator $L$ is surjective on ${\cal C}^{2,\al}_\mu(\Sigma)$ whenever
$\mu > 1$, $\mu \neq 2, 3, \ldots$, or equivalently, whenever
there are no anomalous decaying Jacobi fields and so 
$L$ is injective on ${\cal C}^{2,\al}_{-\mu}(\Sigma)$ for $\mu > 1$. 
\label{def:6.1}
\end{definition}
We cannot preclude the existence of Jacobi fields in ${\cal C}^{2,
\al}_{\mu}(\Sigma)$ for $|\mu| \leq 1$, and in fact these always exist,
at least locally on each end, 
for geometric reasons. They may be exhibited explicitly on the catenoid: 
just as for the Delaunay surfaces, solutions of $L w = 0$ corresponding to the 
eigenvalues of the cross-sectional Laplacian with $|j| \leq 1$ arise from 
translations, rotations and dilations (which substitute for changes in 
Delaunay parameter): 
\begin{proposition}
The Jacobi fields 
\[
\Psi^{0,+}(s) = \tanh s,   \qquad \Psi^{0,-}(s) = s\tanh s - 1,
\]
correspond to vertical translation along the axis of the catenoid, and 
change of the dilation parameter $a$, respectively. The Jacobi fields
\[
\Psi^{1,+}(s,\theta) = \psi^+(s)\cos \theta, \qquad
\Psi^{-1,+}(s,\theta) = \psi^+(s)\sin \theta,
\]
correspond to horizontal translations in the $x_1$ and $x_2$ directions, while
\[
\Psi^{1,-}(s,\theta) = \psi^-(s) \cos \theta, \qquad
\Psi^{-1,-}(s,\theta) = \psi^-(s)\sin \theta,
\]
correspond to rotations about the $x_2$-axis and $x_1$-axis, respectively. Here
\[
\psi^+(s) = \frac{1}{\cosh s}    \qquad \psi^-(s) =  \frac{s}{\cosh s}
+ \sinh s.
\]
\label{pr:6.2}
\end{proposition}
{\bf Proof:} As with the analogous statement in the Delaunay case,
these Jacobi fields may be computed by finding the parametrizations of the 
one-parameter family of minimal surfaces in each case and differentiating
to get the deformation vector field, the inner product of which with
the unit normal of ${\cal C}_a$ yields the appropriate expression.
We leave the details to the reader.
\hfill $\Box$

Jacobi fields asymptotic to these exist on the ends of any $k$-noid. 
Let $\tilde{\Sigma}$ denote some fixed truncation $\Sigma_{\e_0}$
of the $k$-noid $\Sigma$, and let $E_1, \ldots, E_k$ denote the
components of $\Sigma \backslash \tilde{\Sigma}$. These are in
one to one correspondence with the ends of $\Sigma$, and are
minimal surfaces with boundary. 
\begin{proposition} On each end $E_\ell$ of $\Sigma$, there
exists a six-dimensional space of functions $\Psi_\ell^{j,\pm}$, 
$j = 0, \pm 1$, such that each $L \Psi_\ell^{j,\pm} = 0$, and
which are asymptotic to the corresponding model Jacobi fields 
$\Psi^{j,\pm}$ for the catenoid ${\cal C}_{a_\ell}$ modelling $E_\ell$ 
in the sense that
\[
\left| \Psi_\ell^{j,+}(s,\theta) - \Psi^{j,+}(s,\theta)\right|
\le C \, e^{-(j+2)s},
\]
\[
\left| \Psi_\ell^{j,-}(s,\theta) - \Psi^{j,-}(s,\theta)\right|
\le C \, s \, e^{(j-2)s}.
\]
\label{pr:6.3}
\end{proposition}
{\bf Proof:} These new Jacobi fields are produced by the same
geometric process, namely forming the families of minimal 
surfaces with boundary, $E_\ell(\eta)$, by translating, rotating
or dilating $E_\ell$, and then differentiating with respect
to the parameter $\eta$ at $\eta = 0$, and taking the inner
product of the resulting vector field along $E_\ell$ with
the unit normal. The statement about asymptotics is obtained 
from the fact that $E_\ell$ is a normal graph over
${\cal C}_{a_\ell}$ of a function $\phi_\ell$ which
decays like $e^{-2s}$.  \hfill  $\Box$ 

\subsection{Moduli space theory}

Following these preliminaries, we now briefly sketch the moduli space 
theory for $k$-noids. This was developed by Perez and Ros \cite{PR}
at around the same time that the very similar moduli space
theory was set down for solutions of the singular Yamabe problem 
and CMC surfaces in \cite{MPU} and \cite{KMP}, using slightly different 
(but equivalent) methods. The parallels between the three problems are
discussed carefully in \cite{MPo1}. We state the results following these 
latter three papers. For $\Sigma \in {\cal H}_{g,k}$, define the 
$6k$-dimensional space 
\[
W = \oplus_{\ell = 1}^k W_\ell, \qquad
\mbox{where} \qquad W_\ell = \{\eta_\ell \Psi_\ell^{j,\pm}, j = 0, \pm 1\}.
\]
and where $\eta_\ell$ is a cutoff function vanishing on $\tilde\Sigma$ and 
equalling one outside of a slight enlargement of this truncation on the
end $E_\ell$.
At least around non-degenerate $k$-noids, the moduli space theory
is based on the implicit function theorem. For this one requires surjectivity 
of $L$ on some geometrically natural function spaces, but unfortunately the 
spaces on which $L$ is surjective in Proposition~\ref{pr:6.1} have
positive exponential weight, hence are ill-suited for the nonlinear operator. 
To remedy this one uses the following more refined linear result.
\begin{proposition} Suppose $\Sigma$ is nondegenerate. Fix $\mu$ with 
$1 < \mu < 2$. Then the mapping
\begin{equation}
L: {\cal C}^{2,\al}_{\mu}(\Sigma) \oplus W \longrightarrow
{\cal C}^{0,\al}_{\mu-2}(\Sigma),
\label{eq:6.4}
\end{equation}
is surjective. Its nullspace, $B \equiv B_\Sigma$ (which we call the bounded
nullspace) is $3k$-dimensional.
\label{pr:6.4}
\end{proposition}
The proof is essentially identical to the one in \cite{MPU} and \cite{KMP},
although the linear theory here is more elementary than the analysis on 
asymptotically periodic ends in those papers. The dimension
count for $B$ is obtained by a relative index theorem (which is
essentially equivalent to the Riemann-Roch theorem). 

To make sense of the mean curvature operator $N$ on elements of the domain
space in (\ref{eq:6.4}), we use that elements of $W$ correspond to geometric 
motions. Thus $N(u',\tilde{u})$ calculates the mean curvature of the normal
graph of the function $\tilde{u} \in {\cal C}^{2,\al}_{-\mu}$ over the
surface $\Sigma_{u'}$ obtained by slightly deforming the ends of $\Sigma$
in the manner prescribed by the components of $u' \in W$. (More
specifically, one considers an `exponential map' from a neighbourhood
of $0$ in $W$ to a space of surfaces deforming $\Sigma$, such that
the derivatives of the families of surfaces $\Sigma(\lambda u')$, for
$u' \in W$, at $\lambda = 0$ equal $u'$.) Proposition~\ref{pr:6.4} then
states that the differential of this map $N$ is surjective at $(0,0)$
when $\Sigma$ is nondegenerate, and so the first part of the following is 
a trivial application of the standard implicit function theorem: 
\begin{corollary}
In the neighbourhood of any one of its nondegenerate points, the 
moduli space ${\cal H}_{g,k}$ of $k$-noids of genus $g$ is a real 
analytic manifold of dimension $3k$. In the neighbourhood of
an arbitrary point, it has the structure of a locally defined
(possibly singular) real analytic variety. 
\label{cor:6.1}
\end{corollary}
Perez and Ros note that the second part of this result follows from the 
general theory of Weierstrass representations of $k$-noids. It may also 
be proved in a manner more consistent with the first part by the 
Kuranishi method, as in \cite{KMP}, cf. also \cite{MPo1}.

\section{Truncated $k$-noids and their deformations}

In this section we first introduce the compact truncations $\Sigma_\e$ 
which fill out the $k$-noid $\Sigma$ as $\e \rightarrow 0$. 
These are the building blocks occupying the central portion of the 
surfaces we shall construct. The reason for introducing them
is that there are no surfaces of mean curvature one which may be written 
as normal graphs over all of $\Sigma$, but there are many which are 
graphs over any one of the $\Sigma_\e$. Next we study a natural 
boundary problem for the Jacobi operator on these compact surfaces 
and analyze its behaviour as $\e$ tends to zero. Finally, we
introduce a finite dimensional family of deformations of the Jacobi
operator, corresponding to the elements of $W$, and show that
the preceding linear analysis carries over to the operators in this family.

\subsection{The Jacobi operator on truncated $k$-noids} 

We start by defining the truncations $\Sigma_\e$. Recall that
each end $E_\ell$ of $\Sigma$ admits the parametrization
\begin{equation}
\x_\ell (s,\theta)= a_\ell \,( \cosh s \cos \theta + O(e^{-3s}), \cosh s 
\sin \theta + O(e^{-3s}), s+ O(e^{-2s})).
\label{eq:paramav}
\end{equation}
We simply define $\Sigma_\e$ to be the union of the compact piece $K$ 
of $\Sigma$ and the portion of each of the ends up to $s = 
S_\e /8$ (which is of the order $-\frac{1}{4} \log \e$). 

Preliminary to the nonlinear analysis, in the next section,
of the family of surfaces of constant mean curvature one 
which are normal graphs over the $\Se$, we shall
require information about a certain inhomogeneous boundary problem 
for the Jacobi operator $L$ on $\Se$, in particular its solvability 
and the uniformity of this solution with respect to $\e$. 

Before stating this boundary problem precisely, we digress briefly. 
Set 
\[
X(\Se) \equiv {\cal C}^{2,\al}(\Se).
\]
If $u \in X(\Se)$, then we let $(\bC)_\e(u)$ denote its Cauchy data
on $\del \Se$. Thus 
\[
(\bC)_\e(u) \in Z(\dSe) \equiv {\cal C}^{2,\al}(\del \Sigma_\e) 
\oplus {\cal C}^{1,\al}(\del \Sigma_\e).
\]
Since 
\[
{\cal C}^{2,\al}(\Se) = \left. {\cal C}^{2,\al}_{\mu}(\Sigma) \right|_{\Se},
\]
we may use the restriction of $\|\cdot \|_{2,\al,\mu}$ as a norm on $X(\Se)$. 

For later use we also record that this space admits a decomposition 
\[
X(\Se) = W \oplus X(\Se)'',
\]
where, by a slight abuse of notation, $W$ here represents the restrictions 
to $\Se$ of elements of the true `global' deficiency space, and 
\[
X(\Se)'' = \left\{ w \in X(\Se):  \left. w \right|_{E_{\ell}} \equiv w_{\ell}
(s,\theta) \in \mbox{\rm span} \ \left\{\chi_j(\theta)\right\}_{|j| \geq 2} 
\ \mbox{\rm for }\  s_\ell \leq s \leq S_\e /8  \right\}. 
\]
If $X$ is any of the space of functions we consider on all of
$\Sigma$, for example ${\mathcal C}^{2,\al}_{\mu}(\Sigma)$ or
${\mathcal C}^{2,\al}_{-\mu}(\Sigma)$, then we let
$X''$ denote the finite codimension subspace defined in 
the analogous manner (omitting only the restriction
$s \leq S_\e /8 $ in the definition). Most commonly,
$X$ will be ${\mathcal C}^{2,\al}_{-\mu}(\Sigma)$, and so
we write
\[
X_W(\Sigma) = {\mathcal C}^{2,\al}_{-\mu}(\Sigma) \oplus W ,
\]
for brevity. Notice that the bounded nullspace $B$ is
a subspace of $X_W(\Sigma)$. 

Next, if $v \in X(\Se)$, denote its Cauchy data at $\dSe$ by ${\bC}_\e(v)
\in Z(\dSe)$. Similarly, if $v \in X_W(\Sigma)$, 
and its decomposition is written $v = w + v''$, then we regard 
its component $w \in W$ as its Cauchy data at infinity,
and write $w = {\bC}_0(v)$. 

There is a natural (weak) symplectic structure on $Z(\dSe)$ given by 
\[
\omega_\e (\phi, \psi) \equiv \int_{\dSe} \left(\phi_0 \psi_1
- \psi_0 \phi_1 \right)\,d\sigma,
\]
if $\phi = (\phi_0,\phi_1),\ \psi = (\psi_0,\psi_1) \in Z(\dSe)$, 
and where $d\sigma$ is the length form of $\dSe$. (Note that $\omega_\e$ 
does not induce an isomorphism between $Z$ and its dual; fortunately,
this is unimportant for our purposes.) If $u, v\in X(\Se)$ and
$\phi = {\bC}_\e(u),\ \psi = {\bC}_\e(v)$, then by Stokes' theorem,
\begin{equation}
\omega_\e(\phi,\psi) = \int_{\dSe}\left(\frac{\del u}{\del \nu}v - u
\frac{\del v}{\del \nu}\right)\,d\sigma = \int_{\Sigma_\e} (Lu) v - u (Lv).
\label{eq:6.5}
\end{equation}
Here, of course, $\nu$ the appropriately oriented unit normal of $\dSe$.
The expression on the right does not depend on the particular extensions
$u$ and $v$ of $\phi$ and $\psi$. However, fixing $u, v \in
X_W(\Sigma)$, then this expression has a limit as $\e \rightarrow 0$
because now both $Lu$ and $Lv$ decay exponentially on the ends of 
$\Sigma$. In fact, write 
\[
u = u' + u'' \qquad \mbox{where}\qquad u' = \sum_{\ell,j,\pm} 
u_\ell^{j,\pm}  \Psi_\ell^{j,\pm}, \qquad  u'' \in X(\Sigma)'' ,
\]
and similarly $v = v' + v''$; then 
\begin{equation}
\lim_{\e \rightarrow 0} \omega_\e(\phi,\psi) \equiv \omega(\phi,\psi) 
= \ha \sum_{\ell = 1}^k \sum_{j = 0, \pm 1} \left(u_\ell^{j,+}v_\ell^{j,-} - 
u_\ell^{j,-} v_\ell^{j,+}\right). 
\label{eq:6.6}
\end{equation}
Abusing notation in our customary manner, we identify $\omega(\phi,\psi)$ 
with $\omega(u',v')$ and even with $\omega(u,v)$; this is the induced 
symplectic form on $W$, which by this computation is the `standard' 
one with respect to the basis $\{\Psi_\ell^{j,\pm}\}$. 

From (\ref{eq:6.5}), if $u, v \in B$, then $\omega(u,v) = 0$. Since
$\dim B = \ha \dim W$, we conclude that $B$ is Lagrangian in 
the symplectic vector space $(W,\omega)$. From Corollary~\ref{cor:6.1}, 
$B$ is identified with $T_{\Sigma} {\cal H}_{g,k}$, and so we have shown 
that this moduli space inherits at least the infinitesimal structure 
of a Lagrangian submanifold of some larger symplectic manifold. 
This was discussed both in \cite{PR2} and in \cite{MPU}, 
\cite{KMP}. 

We now set up the boundary problem for the Jacobi operator on the 
surfaces $\Sigma_\e$. That the mapping
\[
L: \{u \in {\mathcal C}^{2,\al}(\Se): \left. u \right|_{\dSe} = 0\} 
\longrightarrow {\mathcal C}^{0,\al}(\Se) ,
\]
is surjective when is $\e$ is sufficiently small is fairly 
easy to establish. Unfortunately, the norm of the inverse is not
uniformly bounded as $\e \rightarrow 0$. This happens for
a good reason: the range of the inverse may be too close to
the restriction of the bounded nullspace $B$ to $\Se$. Therefore
we impose a boundary condition, the corresponding solution operator
for which does have the uniformity we need later. 
 
First select a $3k$-dimensional subspace $\tB \subset W$ which is Lagrangian 
with respect to $\left. \omega\right|_W$, and which is transverse to $B$. 
There are many ways to do this, of course, but the choice is irrelevant! 
The space $(B \oplus \tB,\omega)$ is then a $6k$-dimensional dimensional 
symplectic subspace of $Z(\dSe)$, which is a small perturbation 
of $W$ when $\e$ is small. Continuing our practice of splitting into low 
and high eigenspaces, we write $B \oplus \tB$ as $Z(\dSe)'$ because,
up to an error which decreases with $\e$, it corresponds to the span 
in $Z(\dSe)$ of the eigenfunctions $\{\chi_j(\theta)\}_{|j| \leq 1}$. 
For its complement we use 
\[
Z(\dSe)'' \equiv \left[{\cal C}^{2,\al}(\dSe) \oplus {\cal C}^{1,\al}(\dSe)
\right]'',
\]
the span of the eigenfunctions $\{\chi_j(\theta)\}_{|j|>1}$. Thus if
\[
\phi \in Z(\dSe) \quad \mbox{\rm then} \quad \phi = \phi' + \phi'' \quad
\mbox{\rm with }\quad  \phi' \in Z(\dSe)',\quad  \phi'' \in Z(\dSe)'' ,
\]
and these components split further as 
\[
\phi' = \phi_B + \phi_{\tB} \in B \oplus \tB \quad \mbox{\rm and } \quad 
\phi'' = (\phi_0'',\phi_1'') \in {\mathcal C}^{2,\al}(\Se) \oplus
{\mathcal C}^{1,\al}(\Se).
\]
Finally, define the projection
\[
\begin{array}{rlll}
\Pi_{B}: Z(\dSe) & \longrightarrow & B  \\[2mm]
\phi & \longrightarrow  & \phi_{B},
\end{array}
\]
which has nullspace $\tB \oplus\left[ {\mathcal C}^{2,\al}(\dSe) \oplus 
{\mathcal C}^{1,\al}(\dSe)\right]''$. Both $\Pi_{B}$ and the space on 
which it acts depend on $\e$, but we omit this from the notation for 
simplicity. We let $\Pi_{B,0}$ denote the projection sending $\phi \in 
Z(\dSe)$ to $\phi_B + \phi_0''$. 

Finally, for $f \in {\mathcal C}^{0,\al}_{\mu-2}(\Sigma)$ and $\phi_0'' 
\in {\mathcal C}^{2,\al}(\dSe)$, consider the boundary problem
\begin{equation}
\left\{
\begin{array}{rlll}
       Lu & = & \left. f \right|_{\Se}        & \mbox{in} \qquad \Se \\[2mm]
\Pi_{B,0}({\bC}_\e(u)) & = & \phi_0'' & \mbox{on} \qquad \dSe . 
\end{array}
\right.
\label{eq:6.8}
\end{equation}

\begin{proposition} There exists an $\e_0>0$ such that if $\e < \e_0$,
then (\ref{eq:6.8}) has a unique solution $u \in X(\Se)$. Furthermore,
there exists a constant $c>0$, independent of $\e \in (0, \e_0)$, such that 
\[
\|u\|_{2,\al,\mu} \leq c\left( \e^{\mu/4} \|\phi_0''\|_{2,\al} +
\|f\|_{0,\al,\mu-2}\right).
\]
\label{pr:6.5}
\end{proposition}

{\bf Proof:}  To start, use a bounded extension operator in the
${\mathcal C}^{2,\al}(\dSe) \rightarrow {\mathcal C}^{2,\al}_{\mu}
(\Sigma)$ in the usual way to reduce to the case where $\phi_0'' = 0$. 

Our first claim is that when $\e$ is sufficiently small, the only 
solution of the problem (\ref{eq:6.8}) with $f=0$ and $\phi_0'' = 0$
is $u=0$. Granting this for the moment, then because the range of 
$\Pi_B$ is Lagrangian with respect to the forms $\omega_\e$, the 
problem (\ref{eq:6.8}) is self-adjoint, hence the inhomogeneous
problem (\ref{eq:6.8}) has a unique solution. 

Our second claim is that there exists a constant $c >0$, independent 
of $\e >0$ such that 
\begin{equation}
||u||_{2, \alpha, \mu} \leq c ||f||_{2, \alpha, \mu-2}.
\label{eq:6.82}
\end{equation}

The proposition follows from these two claims. As will be seen
from the following argument, the proof of the first claim is
a special case, or at least follows directly from, the proof
of the second claim. Therefore, we concentrate on proving the estimate
(\ref{eq:6.82}), which we do by contradiction.  If it were to fail,
then there would exist a sequence $\e_\ell \rightarrow 0$, functions 
$f_\ell \in {\mathcal C}^{0,\al}_{\mu-2}$ and solutions $u_\ell 
\in X(\Se)$ of (\ref{eq:6.8}) such that
\begin{equation}
\|u_\ell\|_{2,\al,\mu} \equiv 1  \qquad \mbox{and}\qquad  
\|f_\ell\|_{0,\al,\mu-2} \rightarrow 0.
\label{eq:6.88}
\end{equation}

We make a preliminary adjustment of these functions. For each
end $E_m$ of $\Sigma$, find a solution $\bu_\ell \in 
{\mathcal C}^{2,\al}_{\mu}(E_m)$ of 
\[
\left\{
\begin{array}{rllll}
L\bu_\ell & = & \left. f_\ell\right|_{E_m} \\[2mm] 
\left. \bu_\ell \right|_{\del E_m} & = & 0.
\end{array}
\right.
\]
This solution is not unique, but if we choose it orthogonal to 
the nullspace of this problem, then it is well-defined.

We claim that the ${\mathcal C}^{2,\al}_{\mu}$ norm of 
$\bu_\ell$ is bounded. This is essentially just a local form of 
(\ref{eq:6.88}), and may be proved by contradiction as well. If it 
were false, then we could renormalize to make $\sup e^{-\mu s}
|\bu_\ell|=1$ on $E_m$ one. We call the new function $\bu_\ell$ as well. 
Assume that the maximum of $e^{-\mu s}|\bu_\ell|$ occurs
at some point $p_\ell = (s_\ell,\theta_\ell)$.  Then $s_\ell > s_m$,
where $s_m$ gives the cylindrical coordinate of the `inner' boundary
of $E_m$. Define now 
\[
\bw_\ell(s,\theta) = e^{-\mu s_\ell} w_\ell(s+s_\ell,\theta),
\]
which is defined on $(s_m- s_\ell, S_\e /8 - s_\ell) \times S^1$ 
and satisfies $e^{-\mu s}|{\bar{w}}_\ell| \le 1$ on this domain with 
the supremum of $1$ attained at $(0,\theta_\ell)$, 
and of course vanishes at both boundary components of its domain
of definition. 

Pass to a subsequence to obtain a limit $\bw \in {\mathcal C}^{2,
\al}_{\mu}$ which is defined on $(\zeta^-,\zeta^+)\times S^1$ for some 
$\zeta^{\pm} \in {\RR} \, \cup \{\pm \infty\}$, and satisfies Dirichlet 
conditions at either of the boundaries if they remain finite, i.e. if 
$|\zeta^{\pm}|<\infty$. If $s_\ell$ had remained bounded, so $\zeta_- > 
-\infty$, then $L\bw = 0$ on the infinite end $[-\zeta_-,\infty) \times
S^1$ and vanishes at $s=-\zeta_-$, so that $\bw$ would be in the
nullspace of $L$, which is a contradiction. Therefore $s_\ell \rightarrow
\infty$, $\zeta_- = -\infty$, and $(\del_s^2 + \del_\theta^2)\bw = 0$. But
solutions of this equation are sums of exponentials, multiplied by 
eigenfunctions on the cross-section, and one checks readily that no 
there is no solution such that $|\bw|\leq e^{\mu s}$ either on all 
of ${\RR}$ or on $(-\infty,\zeta^{+}]$. This gives a contradiction
again and our claim is proved.

By subtracting off $\chi \bu_\ell$ from $u_\ell$, where $\chi$ is
smooth, vanishes for $s \leq s_m$  and equals one for $s \geq s_m + 1$
on each $E_m$, we have now reduced to the case where $f_\ell$
is compactly supported. We do not change the names of these
adjusted functions, and still assume that they satisfy (\ref{eq:6.88}).

Now we continue in a very similar fashion as before. Select a smooth 
positive function $d$ on $\Sigma$ which agrees on each end with the 
coordinate function $s$, so that 
\[
\sup e^{-\mu d}|u_\ell| = 1.
\]
This supremum is attained at some point $p_\ell \in \Sigma_{\e_\ell}$. 

If $p_\ell$ were to tend to infinity, then we could assume that it did 
so within one end $E_m$. We then apply an argument identical to 
the one just above to reach a contradiction. 

The final case is when $p_\ell$ remains in some fixed compact set in 
$\Sigma$, which means that we can extract a subsequence converging to some 
function $u \in {\cal C}^{2,\al}_{\mu}(\Sigma)$ such that $Lu = 0$ and
$u \not\equiv 0$. Hence $u$ is an element of the bounded nullspace $B$
(and so in particular, lies in $W \oplus {\mathcal C}^{2,\al}_{-\mu}
(\Sigma)$).  Let ${\bC}_0(u) = \phi \in B \subset W$.  Then, choose
any $\psi \in \tB$ and an extension $v \in W + {\mathcal C}^{2,\al}_{-\mu}
(\Sigma)$ of it to $\Sigma$ so that $\psi = {\bC}_0(v)$, by definition, 
and $Lv$ is compactly supported.

Then we compute that 
\[
\begin{array}{rlll}
\omega(\phi,\psi) & = &   \displaystyle \int_\Sigma \left( (Lu) v - u(Lv) \right)\\[2mm]
 \, & = &   \displaystyle  \lim_{\ell \rightarrow \infty} \int_{\Sigma_{\e_\ell}}
\left((Lu_\ell) v - u_\ell (Lv)\right) 
= \lim_{\ell \rightarrow \infty} \omega_\e({\bf C}_{\e_\ell}(u_\ell), 
{\bf C}_{\e_\ell}(v) = 0.
\end{array}
\]
The first equality is obvious, while for the second one, the compact
supports of $Lu_\ell = f_\ell$ and $Lv$ ensure the convergence. The
third is again obvious, while the final equality holds because 
$\tB$ is Lagrangian, so the contribution from $Z'$ is zero, while 
the contribution from $Z''$ tends to zero with $\e$. But $\phi \in B$ 
while $\psi$ is an arbitrary element of $\tB$, and the symplectic 
pairing between $B$ and $\tB$ is nondegenerate. This proves that $\phi$ 
and hence $u$ vanishes identically, which is a contradiction.

The proof is complete in all cases. \hfill $\Box$

\subsection{Deformations of $\Se$}

Now we shall take up the task of defining slightly different truncations
of the scaled surface $\e \, \Sigma$, which we shall call $
\tilde{\Sigma}_{\e, {\cal P}}$ (the ${\cal P}$ here refers
to a parameter set which we shall define below), which will be more
convenient later. In the next subsection we shall also consider the 
Jacobi operators which correspond to writing nearby surfaces
as graphs using vector fields which are small deformations of the
normal vector field on $\tilde{\Sigma}_{\e, {\cal P}}$. 

Fix one end $E_\ell$ of $\Sigma$, and recall its parametrization
(\ref{eq:paramav}). The end $\e E_\ell$ can therefore be paramererized as
\begin{equation}
\x_\ell (s,\theta)= \e_\ell \,( \cosh s \cos \theta + O(e^{-3s}), \cosh s 
\sin \theta + O(e^{-3s}), s+ O(e^{-2s})).
\label{eq:paramavbis}
\end{equation}
where here and later, we use the notation
\[
\e_\ell \equiv a_\ell \, \e.
\]
There is a Delaunay surface ${\cal D}_\ell$, 
with Delaunay parameter $\e_\ell$, which `best fits' the model catenoid
for $\e E_\ell$ near the region where $s = S_{\e_\ell}/8$. 
It has the parametrization
\[
\x_{{\cal D}_\ell}(s,\theta)=(\tau_\ell e^{\sigma_\ell(s)} \cos\theta,
\tau_\ell e^{\sigma_\ell(s)} \sin \theta, k_\ell(s)),
\]
and has unit normal 
\[
\nu_{{\cal D}_\ell}(s,\theta) = (- \tau_\ell \cosh( {\sigma_\ell(s)}) 
\cos \theta, - \tau_\ell \cosh( {\sigma_\ell(s)}) \sin \theta, 
(\sigma_{\ell})_s(s)).
\]

The deformations of ${\cal D}_\ell$ are parametrized by translations 
orthogonal to the axis, translations along this axis, rotations of this 
axis and finally changes in the Delaunay parameter. We label these by 
\[
{\cal P}^\ell = (t_1^\ell,t_2^\ell , r_1^\ell, r_2^\ell, d_\ell , \delta^\ell),
\]
respectively. All these parameters lie in some small neighbourhood of 
zero. (The $r$'s are identified with some small neighbourhood of
the identity in the space of rotations fixing the $x_3$ axis; the
exact manner is not important, but to be definite, we suppose
that the diffeomorphism is given by the exponential map in $SO_3$ 
orthogonal to the copy of $SO_2$ which is the stabilizer of
that axis, followed by the projection to $SO_3/SO_2$.)
The full parameter set for all ends of $\e \, \Sigma_\e$ is
\[
{\cal P} = ({\cal P}^1, \ldots,{\cal P}^k),
\]
We also set
\[
\tilde{t}= (t^1_1, \ldots, t^k_2) \in {\RR}^{2k}, \qquad \tilde{r} = 
(r^1_1, \ldots, r^k_2) \in {\RR}^{2k},
\]
\[
\tilde{d}= (d^1, \ldots, d^k) \in {\RR}^{k}, \qquad\mbox{and} \qquad
\tilde{\delta}= (\delta^1, \ldots, \delta^k) \in {\RR}^{k},
\]
and the rigid motion determined by $(t_1^\ell,t_2^\ell,r_1^\ell,r_2^\ell,
d^\ell)$ will be denoted ${\cal R}_{(t^\ell,r^\ell,d^\ell)}$. The norm on 
these parameter sets which arises naturally below is given by 
\[
\|{\cal P}\| = \| (\tilde{t}, \tilde{r}, \tilde{d}, \tilde{\delta})\| 
\equiv  \e^{1/4} \| \tilde{t}\| + \e^{3/4}\|\tilde{r}\| + \|\tilde{d}\| +
(\log \frac{1}{\e})  \, \| \tilde{\delta} \|.
\]

The Delaunay surface associated to the set of (small) deformation parameters 
${\cal P}_\ell$ will be denoted ${\cal D}_{{\cal P}^\ell}$, and its induced
parametrization and unit normal will be called $\x_{{\cal P}^\ell}$ and 
$\nu_{{\cal P}^\ell}$, respectively. This surface has Delaunay parameter 
$\e_\ell + \delta^\ell$. 

We come now to the main point, which is to write a neighbourhood of 
$\e \, ( E_\ell \cap \dSe ) $ as a normal graph over each 
${\cal D}_{{\cal P}^\ell}$, and to obtain estimates on the graph function.

\begin{proposition}
Fix $\kappa \in (1, \frac{3}{2})$. Then, for all parameter sets ${\cal P}$ with $\|{\cal P}\| \leq 
\e^\kappa $, there is a diffeomorphism 
$\Psi(s,\theta) = (s',\theta')$ from $(-2 + S_{\e_\ell} /8, 
S_{\e_\ell}/8)\times S^1$ onto its image, satisfying
\[
\|\Psi(s,\theta) - (s,\theta)\| = O(\e^{\kappa - 1}),
\]
and we have 
\begin{equation}
\x_\ell (s,\theta) = \x_{{\cal P}^\ell}(s',\theta') + \hat{w}_0(s', \theta')
\nu_{{\cal P}^\ell}(s',\theta'),
\label{eq:par}
\end{equation}
for all $\e < \e_0$, where $\e_0$ depends only on $\kappa$. 
The graph function $\hat{w}_0$ here is of the form
\[
\begin{array}{rlll}
\hat{w}_0 (s', \theta') & = & -\displaystyle \frac{1}{\cosh s'} ( t_1^\ell 
\cos \theta' + t_2^\ell \sin \theta' )  -  ( r_1^{\ell} \cos \theta' + 
r^{\ell}_2 \sin \theta' ) \, \e_\ell \,  \cosh s'  \\[2mm]
& \, & + d^{\ell}  +\delta^\ell s'  +O(\e^{3/2} + \e^{2\kappa -1}).
\end{array}
\]
\label{pr:modif}
\end{proposition}

In other words, we are writing a neighbourhood of $\e \, ( E_\ell \cap \dSe) $ as a normal
graph over each of the family of nearby model Delaunay surfaces, up to the 
reparametrization given by the diffeomorphism $\Psi$. 

{\bf Proof :}
This follows from a computation similar to the one we have already 
done in the proof of Proposition~\ref{pr:4.1}. Recall that in the
range $ s\in [- 4 + S_\e /8, 4 + S_\e/8]$, we have the 
expansions
\[
k(s) = \e s + O(\e^{3/2}), \qquad \qquad  \tau e^{\sigma (s)} = \e \cosh s 
+ O( \e^{5/4}),
\]
and 
\[
\tau \cosh \sigma(s) = \frac{1}{\cosh s} + O(\e^{3/4}), \qquad \qquad \del_s 
\sigma (s) = 1 +O(\e^{1/2}),
\]
which follow from (\ref{eq:colexp}). Here and below $O(\e^{\gamma})$ 
will denote functions of $(s,\theta)$ all derivatives of which are 
bounded by constant multiples of $\e^\gamma$. 

It will be most convenient to apply the transformation 
${\cal R}^{-1}_{(t^\ell,r^\ell,d^\ell)}$ to both sides of (\ref{eq:par}).
On the one hand, from (\ref{eq:paramav}), the parametrization for
${\cal R}^{-1}_{(t^\ell,r^\ell)}(\x_\ell (s,\theta))$ is given by
\[
\begin{array}{rllll}
(s,\theta)\longrightarrow & \left(\e_\ell\cosh s \cos\theta -t^\ell_1
+ O(\e^{\kappa + 1/4}\log \e), \right. \\[2mm]
\, & \   \e_\ell \cosh s \sin  \theta -t^\ell_2 + O(\e^{\kappa + 1/4}
\log \e), \\[2mm]
\, &  \left. ( r_1^{\ell} \cos \theta + r^{\ell}_2 \sin \theta ) \e_\ell 
\cosh s  + \e_\ell s  - d^{\ell} +O(\e^{3/2}) \right)
\end{array}
\]
for $s$ in this range.

On the other hand, ${\cal R}^{-1}_{(t^\ell,r^\ell,d^\ell)}
(\x_{{\cal P}^\ell}(s',\theta') + \hat{w}_0 (s', \theta')
\nu_{{\cal P}^\ell}(s',\theta'))$ is parameterized by
\[
\begin{array}{rllll}
(s', \theta') & \longrightarrow  \\[2mm]
&\big(((\e_\ell+\delta^\ell) \cosh s' - \displaystyle \frac{1}{\cosh s'} \hat{w}_0 
(s',\theta'))\cos\theta' + O(\e^{5/4}) + O(\e^{3/4}) \hat{w}_0(s',\theta'),
\\[2mm]
\, & \,((\e_\ell+\delta^\ell) \cosh s' - \displaystyle \frac{1}{\cosh s'} \hat{w}_0 
(s', \theta') ) \sin \theta' + O(\e^{5/4}) + O(\e^{3/4}) \hat{w}_0(s',\theta'),
\\[2mm]
\, &\,(\e_\ell+\delta^\ell) s' + \hat{w}_0(s', \theta')+  O(\e^{3/2}) 
+ O(\e^{1/2}) \hat{w}_0(s', \theta') \big),
\end{array}
\]
again for $s'$ in this range.

Equating the third coordinates, we already find that $\hat{w}_0(s',\theta') 
= \e(s-s') + O(\e^\kappa)$. Assuming that $|s'-s|$ is at least
bounded, this gives $\hat{w}_0(s',\theta') = O(\e)$. Similar estimates 
hold for its derivatives. Now, writing out the equality of the three 
coordinates in turn gives
\[
\begin{array}{rll}
\e_\ell \cosh s \cos \theta   &  = & t^\ell_1 + ((\e_\ell+\delta^\ell)
\cosh s' - \displaystyle \frac{1}{\cosh s'} \hat{w}_0 (s', \theta') ) \cos
\theta' + O(\e^{5/4}) \\[2mm]
\e_\ell \cosh s \sin  \theta  &  = & t^\ell_2 + ((\e_\ell+\delta^\ell) 
\cosh s' - \displaystyle \frac{1}{\cosh s'} \hat{w}_0 (s', \theta') ) \sin 
\theta' + O(\e^{5/4}),
\end{array}
\] 
and 
\[
(r_1^{\ell} \cos \theta + r^{\ell}_2 \sin \theta ) \e_\ell \cosh s  + 
\e_\ell s  - d^{\ell}  =  (\e_\ell+\delta^\ell) s' + \hat{w}_0(s', \theta') 
+O(\e^{3/2} \log \e) .
\]
Using the preliminary estimate on $\hat{w}_0$, we conclude that 
\[
|s-s'|\leq c \, \e^{\kappa-1} , \qquad |\theta' -\theta|\leq c\,  \e^{\kappa-1},
\]
and then, reinserting this information back into the third equality, that
\[
|\hat{w}_0|\leq c \, \e^{\kappa},
\]
along with its derivatives. The third identity gives
\[ 
\begin{array}{rlll}
\hat{w}_0(s',\theta') & = & -\e_\ell \,  (s-s') -  ( r_1^{\ell}  \cos \theta' + 
r^{\ell}_2  \sin \theta' ) \, \e_\ell  \cosh s' \\[2mm]  
&\, & + d^{\ell}+\delta^\ell s' +O(\e^{3/2} + \e^{2\kappa -1}),
\end{array}
\]
while from the first two identitites we get
\[
\e_\ell (\cosh s - \cosh s' ) =  t_1^\ell \cos \theta + t_2^\ell \sin \theta 
+ O( \e^{5/4} + \e^{2 \kappa -5/4}).
\]
This leads finally to
\[ 
\begin{array}{rlll}
\hat{w}_0(s', \theta') & = & -\displaystyle \frac{1}{\cosh s'} ( t_1^\ell 
\cos \theta' + t_2^\ell \sin \theta' )  -  ( r_1^{\ell} \cos \theta' + 
r^{\ell}_2 \sin \theta' ) \, \e_\ell  \, \cosh s'  \\[2mm]
& \, & + d^{\ell}  +\delta^\ell s'  +O(\e^{3/2} + \e^{2\kappa -1}),
\end{array}
\]
which is the desired expansion. \hfill $\Box$

We may now define the deformation $\tilde{\Sigma}_{\e, {\cal P}}$ 
when the set of deformation
parameters ${\cal P}$ satisfies $\|{\cal P}\| \leq \e^{\kappa}$.
Choose $\e_0$ sufficiently small that $S_{\e_\ell}/8 - 2 > 
s_\ell$ for each $\ell$ whenever $\e < \e_0$. Then for any such $\e$, 
define $\tilde{\Sigma}_{\e, {\cal P}}$ as the union of the central compact portion of $\e \, \Se$ and 
the portion of each end $\e E_\ell$ for 
$s_\ell \leq s \leq - 1 + S_{\e_\ell} /8$ and by the graph of 
\[
(s',\theta') \longrightarrow  \x_{{\cal P}^\ell}(s',\theta') + 
\hat{w}_0 (s',\theta') \, \nu_{{\cal P}^\ell}(s',\theta'),
\]
for $-2 + S_{\e_\ell} /8 \leq s \leq  S_{\e_\ell}/8 $.
\begin{remark}
These definitions are compatible in the  region of overlap, and all
we have done is to slightly alter the boundary of $\tSeP$ so that
it conforms better to the coordinates $(s',\theta')$. 
\end{remark}

\subsection{Deformed Jacobi operators}

For any small parameter set ${\cal P}$, we define on the surface $\tSeP$ 
a vector field $\tilde{\nu}$ which is the unit 
normal vector field away from the boundary, and which is a perturbation of this 
unit normal near to the boundary. More specifically, write
$\e \, E_\ell$ as the graph 
\[
(s,\theta) \longrightarrow \x_{{\cal P}^\ell}(s,\theta) + \hat{w}_0
(s,\theta) \, \nu_{{\cal P}^\ell}(s,\theta),
\]
for all $s \in [-2 + S_{\e_\ell} /8, S_{\e_\ell} /8]$.  Let $\eta(s)$ be a 
smooth cutoff function equal to $1$ for $s\leq - 3/2$ and vanishing for 
$s \geq -1$. Then, for all $ s \in [-2 + S_\e /8, S_\e /8]$, the vector 
$\tilde{\nu}(s,\theta)$ is defined to be the unit normal to the surface 
parameterized by 
\[
(s,\theta) \longrightarrow   \,{\cal R}_{(t^\ell,r^\ell)}
( \x_{{\cal D}_\ell} (s,\theta) + \eta (s +S_{\e_\ell} /8) 
\, \hat{w}_0 (s, \theta) \, \nu_{{\cal D}_\ell} (s,\theta) ),
\]
As desired, $\tilde{\nu}$ is still the unit normal to $\tSeP$ when 
$s \leq -3/2 -S_{\e_\ell} /8$, and equals $\nu_{{\cal P}^\ell}(s,\theta)$ when 
$s \in  [-1 + S_{\e_\ell} /8, S_{\e_\ell} /8]$. 

Any surface near to $\tSeP$ may be parameterized by 
\[
\tSeP \ni p  \longrightarrow  p + w (p) \tilde{\nu}(p),
\]
for some scalar valued function $w$. We need to consider the
equation which $w$ must satisfy in order for this surface
to have constant mean curvature one, which we shall do in
a slightly more general context.

Let $S$ be a regular orientable surface, with unit normal vector
field $\nu$. Suppose that $\bar{\nu}$ is another unit vector field
along $S$ which is nowhere tangential. By the inverse function
theorem, for any $p_0 \in S$ there are neighbourhoods
${\cal U}$ and ${\cal V}$ near $(p_0,0)$ in $S \times \RR$ and 
a diffeomorphism $(\phi(p,s),\psi(p,s))$ from ${\cal U}$ to ${\cal V}$
such that
\begin{equation}
p + s \nu(p) = \phi(p,s) + \psi(p,s)\bar{\nu}(\phi(p,s)).
\label{eq:16.1}
\end{equation}
Here $\phi(p,0) = p$ and $\psi(p,0) = 0$. To determine the
first order Taylor series of these functions in $s$, differentiate
(\ref{eq:16.1}) with respect to $s$ and set $s=0$. This gives
\[
\nu(p) = \frac{\del \phi}{\del s}(p,0) + \frac{\del \psi}{\del s}(p,0)\bar{\nu}(p),
\]
and so, taking the normal component of this, we get
\[
1 = \frac{\del \psi}{\del s}(p,0)\,\nu(p)\cdot \bar{\nu}(p), \qquad
\mbox{\rm or } \quad \frac{\del \psi}{\del s}(p,0) = 1/\big(
\nu(p)\cdot \bar{\nu}(p)\big).
\]
Hence
\[
\psi(p,s) = \frac{s}{\nu(p)\cdot \bar{\nu}(p)} + O(s^2).
\]
On the other hand, taking the tangential component and using this expansion
of $\psi$ yields
\[
0 = \frac{\del \phi}{\del s}(p,0) + \frac{s}{\nu(p)\cdot \bar{\nu}(p)}
\bar{\nu}_t(p),
\]
where $\bar{\nu}_t(p)$ is the tangential component of $\bar{\nu}$. Thus
\[
\phi(p,s) = p - \frac{s}{\nu(p)\cdot \bar{\nu}(p)} \bar{\nu}_t(p) + O(s^2).
\]

Next, any surface which is ${\cal C}^2$ close to $S$ can be parameterized 
either as a normal graph of some function $w$ over $S$, using 
the vector field $\nu$, or as a graph of a different function $\bar{w}$ 
using the vector field $\bar{\nu}$. These functions are related
by
\[
p + w(p) \, \nu(p) = \bar{p} + \bar{w}(\bar{p}) \, \bar{\nu}(\bar{p})
= \phi(p,w(p)) + \psi(p,w(p)) \, \bar{\nu}(\phi(p,w(p)).
\]
Using the expansions above, we see that $\bar{w}(p) = w(p)/
(\nu(p) \cdot \bar{\nu}(p)) + O(\|w \|^2)$. 

The mean curvature operators on these two functions, which we call 
$H_{\nu,w}$ and $H_{\bar{\nu},\bar{w}}$, respectively, are related by 
\begin{equation}
H_{\bar{\nu},\bar{w}}(\bar{p}) = H_{N, w}(p).
\label{eq:16.2}
\end{equation}
Differentiating this with respect to $\bar{w}$ and setting $\bar{w} = 0$, 
we get
\begin{equation}
D_{\bar{w}} H_{\bar{\nu}, 0}(u) = D_{w}H_{\nu,0}((\bar{\nu}\cdot \nu) \, u) + 
\left(\nabla H_{\nu,0} \cdot \bar{\nu}_t\right) \, u,
\label{eq:16.3}
\end{equation}
for any scalar function $u$. In the special case where the surface $S$ has 
constant mean curvature, this reduces to 
\begin{equation}
\bar{L} u \equiv D_{\bar{w}} H_{\bar{\nu},0} (u) = D_{w} H_{\nu,0}( 
(\bar{\nu} \cdot \nu)\, u ) \equiv L( (\bar{\nu} \cdot \nu) \, u).
\label{eq:16.4}
\end{equation}

We apply the previous computation to the present situation.
Denote by $\tLeP$ the linearized mean curvature operator about 
$\tSeP$. Away from $\del \tSeP$, we have 
\[
\tLeP =\displaystyle \frac{1}{\e^2} L,
\]
where $L = \Delta_\Sigma + |A_\Sigma|^2$ is the operator we
have studied in detail. Near $\del \tSeP$ the structure of $\tLeP$ 
is described by the next result, the proof of which follows from the 
expansions given in  Proposition \ref{pr:modif}.
\begin{lemma}
In $\e \, E_\ell$, we can write
\[
\tLeP = \displaystyle \frac{1}{\e^2} L +\hat{L}_{\e,{\cal P}},
\]
where $\hat{L}_{\e,{\cal P}}$ is a second order linear differential operator whose 
coefficients are supported in $[-2 +S_{\e_\ell} /8, S_{\e_\ell} /8] \times 
S^1$ and are bounded by $ \displaystyle \frac{1}{\e^2 e^{2s}} \e^{\kappa-1}$.
\label{le:dlo1}
\end{lemma}

Also, following from the same ideas as in \S 7.1 is the simpler 
\begin{lemma}
In $\e \, E_\ell$, the difference
\[
\displaystyle \frac{1}{\e^2} L - \displaystyle \frac{1}{\e_\ell^2 \cosh^2 s} 
\left(\del_s^2 +\del_\theta^2 \right)
\]
is a second order linear differential operator, the coefficients of
which are bounded by a constant times $\e^{-2} e^{-4s}$ in $[s_\ell , 
S_{\e_\ell} /8]$.
\label{le:dlo2}
\end{lemma}

The proofs of both of these results are left to the reader. 

From these lemmas, we can immediately generalize Proposition \ref{pr:4.3} 
to the deformed Jacobi operators on the surfaces $\tSeP$.
\begin{proposition} 
Fix $\mu$ with $1 < \mu < 2$. Then there exists an $\e_0 >0$, depending only 
on $\mu$, such that whenever $0<\e<\e_0$, there exists a unique solution 
$w \in {\mathcal C}^{2,\al}_{\mu}(\tSeP)$ of the problem 
\begin{equation} 
\left\{ 
\begin{array}{rlll} 
\tLeP w & = & \displaystyle \frac{1}{\e^2} f \quad  & \mbox{in}
\quad \tSeP \\[2mm] 
\Pi'' w & =  &  \phi'' \quad & \mbox{on}\quad  \del \tSeP,
\end{array} 
\right. 
\label{eq:8.1} 
\end{equation} 
for $f \in {\cal C}^{0,\al}_{\mu-2}(\tilde{\Sigma}_\e)$ and $\phi''=
(\phi''_1, \ldots , \phi''_k) \in \Pi''\left( {\cal C}^{2,\al}(S^1)\right)^k$.
The Green and Poisson operators will be denoted $\tilde{G}_{\e,{\cal P}}$
$\tilde{P}_{\e,{\cal P}}$, respectively. The linear maps
\[
\begin{array}{rlll}
\tilde{G}_{\e,{\cal P}}:{\cal C}^{0,\al}_{\mu-2} (\tSeP) 
& \longrightarrow & {\cal C}^{2,\al}_{\mu}(\tSeP),\\[2mm]
\e^{-\mu /4} \tilde{P}_{\e}: \Pi''\left( {\cal C}^{2,\al}(S^1)\right)^k
 & \longrightarrow & {\cal C}^{2,\al}_{\mu}(\tSeP),
\end{array}
\]
are bounded uniformly as $\e \rightarrow 0$.
\label{pr:8.1} 
\end{proposition} 

Following the results of section \S 4.3, we also prove 

\begin{corollary}
Fix $1< \mu < 2$. Then there exists a constant $c>0$ and 
an $\e_0 >0$, depending only on $\mu$, such that for $0<\e<\e_0$, we have
\[
||(\tilde{P}_{\e} - \tilde{P}_{0})(\phi'')||_{2,\al,\mu}
\leq c \, \e^{\mu /4} \left(\e^{(2-\mu)/4} + \e^{\kappa-1}\right)
||\phi||_{2,\al}.
\]
Here, if $\phi''= ( \phi_1, \ldots , \phi_k)  \in \Pi'' \left({\cal C}^{2,\al}
(S^1)\right)^k$, the function  $\tilde{P}_{0}(\phi'') = \tilde{w}_0 $ is 
defined to be equal to $\eta (s-s_\ell) \tilde{w}_\ell$ on each end $\e \, 
E_\ell$ and $0$ elsewhere, where $\eta$ is some cutoff function equal to $0$ 
for $s <0$ and equal to $1$ for $s >1$ and where  $\tilde{w}_\ell$ is the 
unique  solution, in $\left( {\cal C}^{2,\al}_{2} ((-\infty, S_{\e_\ell} /8] 
\times S^1)\right)^k$, of 
the problem
\begin{equation} 
\left\{ 
\begin{array}{rlll} 
\Delta \tilde{w}_\ell & = & 0 \quad &\mbox{in}\quad (-\infty, S_{\e_\ell}/8) 
\times S^1\\[2mm] 
\tilde{w}_\ell  & = & \phi''_\ell\quad &\mbox{on}\quad \{S_{\e_\ell}/8\} 
\times S^1.
\end{array} 
\right. 
\label{eq:8.3}
\end{equation} 
\label{cor:8.2}
\end{corollary}
{\bf Proof :} We start by solving, for each $\ell$, 
\[
\left\{
\begin{array}{llll}
\Delta \tilde{w}_\ell & = & 0  & \mbox{ in} \qquad (-\infty ,S_{\e_\ell}/8) 
\\[2mm]
\tilde{w}_\ell & = & \phi''_\ell & \mbox{ on} \qquad \{S_{\e_\ell}/8\} 
\times S^1 
\end{array}
\right.
\]
There is a unique solution of this equation, which is in ${\cal C}^{2, 
\alpha}_{2} ((-\infty , S_{\e_\ell} /8] \times S^1 )$, and satisfies
\[
||\tilde{w}_\ell ||_{2,\al,2}\leq c \,\e^{1/2}||\phi''||_{2,\al}.
\]
Now truncate these solutions at $s=s_\ell$; this allows one to define 
$\tilde{w}_0$ globally  on $\tSeP$ by setting it equal to
$0$ elsewhere. From Lemmas \ref{le:dlo1} and \ref{le:dlo2} if follows
that on each end $\e\, E_\ell$, the difference 
\[
\tLeP - \displaystyle \frac{1}{\e_{\ell}^2 \cosh^2 s} \left( \del_s^2 +\del_\theta^2 
\right),
\]
is a second order linear differential operator whose coefficients are 
sums of terms which are either bounded by a constant times 
$e^{-2s}$ or are supported in $[-2+S_{\e_\ell}/8, S_{\e_\ell} /8]$ 
and bounded by a constant times $\e^{\kappa}$. Using this, we see that
\[
||\tLeP \tilde{w}_0||_{0,\al,\mu} \leq c (\e^{1/2} +  \e^{\kappa 
 -1 + \mu/4})  \|\phi''\|_{2,\al} .
\]
The result then follows from Proposition~\ref{pr:8.1}.
\hfill $\Box$

\section{CMC surfaces near to the truncated $k$-noids}

Just as we already did for Delaunay surfaces, we would like to
analyze the family of surfaces which are close to each $\tSeP$ 
and which have constant mean curvature $1$. To this end, as in 
(\ref{eq:5.2}), we expand the mean curvature operator to see that, 
for any $\phi \in {\cal C}^{2,\al}(\del \tSeP)$, 
our problem reduces to solve the following boundary value problem
\begin{equation}
\left\{ \begin{array}{rlll}
\tLeP w & = & 1 + \tilde{{\cal Q}} (w) & \qquad \mbox{in} \qquad 
\tSeP \\[2mm]
\Pi'' (w) & = &  \phi'' & \qquad \mbox{on} \qquad  \del \tSeP.
\end{array}
\right.
\label{eq:8.4}
\end{equation}
Here 
\begin{equation}
\tilde{{\cal Q}} (w) = \frac{1}{\e e^s} \tilde{Q} \left(\frac{w}{\e e^s} , 
\frac{\nabla w}{\e e^s}, \frac{ \nabla^2 w}{\e e^s}\right),
\label{eq:8.5}
\end{equation}
in each end $\tilde{E}_\e$, collects all the terms of order
higher than one in $w$. The function $\tilde{\cal Q}$ has partial derivatives 
which are uniformly bounded. Denote by $\tilde{w}_\e$ the solution of 
\[
\left\{ \begin{array}{rlll}
\tLeP \tilde{w}_\e &  =  & 0 & \qquad \mbox{in} \qquad \tSeP \\[2mm]
\Pi'' \tilde{w}_\e &  =  & \phi'' & \qquad \mbox{on} \qquad \del \tSeP,
\end{array}
\right.
\]
which is given by Proposition~\ref{pr:8.1}. By the same Proposition, we can 
also solve 
\[
\left\{ \begin{array}{rlll}
\tLeP \tilde{w}_1 & = &  1  & \qquad \mbox{in} 
\qquad  \tSeP \\[2mm]
\Pi'' \tilde{w}_1 & = & 0 & \qquad \mbox{on} \qquad \del \tSeP .
\end{array}
\right.
\]
We find that
\[
||\tilde{w}_1 ||_{2, \alpha,\mu} \leq c  \e^{3/2 + \mu /4}.
\]

Setting $w = \tilde{w}_\e + \tilde{w}_1 + v$, then it remains to solve 
\[
\left\{ \begin{array}{rlll}
\tLeP v & = &   \tilde{\cal Q} (\tilde{w}_\e + 
\tilde{w}_1 + v) & \qquad \mbox{in} \qquad  \tSeP \\[2mm]
\Pi'' v & = & 0 & \qquad \mbox{on} \qquad \del \tSeP .
\end{array}
\right.
\]
It is sufficient to find a fixed point of the mapping
\begin{equation}
\tilde{\cal K}(v) \equiv \tilde{G}_{\e,{\cal P}}\tilde{\cal Q}(\tilde{w}_\e + 
\tilde{w}_1 + v ),
\label{eq:tfpdel}
\end{equation}
when $\e$ is sufficiently small.
\begin{proposition} 
There exists a constant $c_0 >0$ such that if $||\phi''||_{2,\al} 
\leq c_0\,\e^{3/4}$, then 
\[
||\tilde{G}_{\e}({\tilde{\cal Q}}( \tilde{w}_\e +\tilde{w}_1 ))||_{2,\al,\mu}
\leq c  \e^{-3/4} \left( ||\phi''||^2_{2,\al}+\e^3 \right) \, \e^{\mu /4},
\]
and
\[
||\tilde{G}_{\e}({\tilde{\cal Q}}(\tilde{w}_\e +\tilde{w}_1+ v_2)
-{\tilde{\cal Q}} (\tilde{w}_\e +\tilde{w}_1 +v_1))||_{2,\al ,\mu} 
\leq \ha ||v_2-v_1||_{2,\al,\mu},
\]
for all $v_1, v_2$ in $\tilde{B}_{c_0} \equiv \{v: ||v||_{2,\al,\mu}\leq 
c_0\,\e^{(3+\mu)/4}\}$. Thus, $\tilde{\mathcal K}$ is a contraction mapping
on the ball $\tilde{B}_{c_0}$ into itself, and therefore 
has a unique fixed point $v$ in this ball.
\label{pr:tfixptdel}
\end{proposition}
{\bf Proof :}
We use that 
\begin{equation}
||\tilde{w}_\e||_{2,\al,[s,s+1]} \leq c\, e^{\mu (s-S_{\e_\ell}/8)} \,
||\phi''||_{2,\al} \leq c\,c_0\,\e^{3/4}\,  e^{\mu(s -S_{\e_\ell} /8)},
\label{eq:test1}
\end{equation}
and also that
\begin{equation}
||\tilde{w}_1||_{2,\al,[s,s+1]} \leq c \e^{3/2} e^{\mu(s - S_{\e_\ell} /8)}.
\label{eq:test2}
\end{equation}
These estimates imply that on the end $\e E_\ell$, for $s \in [s_\ell, 
S_{\e_\ell} /8]$, we have 
\[
\displaystyle e^{-\mu s} \left\| \frac{1}{\e e^{s}} \left(\frac{\tilde{w}_\e 
+\tilde{w}_0}{\e e^{s}} \right)^2 \right\|_{2,\al,[s,s+1]} 
\leq c\,  \e^{-3/4}\left( ||\phi''||^2_{2,\al}+\e^3 \right) \, \e^{\mu /4},
\]
and then, from a Taylor expansion we get
\[
e^{-\mu s}|| \tilde{\cal Q}(\tilde{w}_\e +\tilde{w}_1)||_{0,\al,[s,s+1]}
\leq c\,  \e^{-3/4}\left( ||\phi''||^2_{2,\al}+\e^3 \right) \, \e^{\mu /4}
\leq  c\, c_0^2\,\e^{(3-\mu)/4}. 
\]
On the other hand, on the compact piece, we simply have 
\[
e^{-\mu s}|| \tilde{\cal Q}(\tilde{w}_\e +\tilde{w}_1)||_{0,\al,[s,s+1]}
\leq c\,  \e^{-3/4} ||\phi''||^2_{2,\al} \, \e^{\mu /4}
\leq  c\, c_0^2\,\e^{(3-\mu)/4}. 
\]
The other estimate follows in the same way, and the proof is
complete. \hfill $\Box$

As in section \S 4.3 we finally obtain
\begin{corollary}
There exists a constant $c_0 >0$ and an $\e_0>0$ such that, for all $\e
\in (0, \e_0)$ and  for any $\phi'' \in \Pi''\left( 
{\cal C}^{2,\al}(S^1)\right)^k $ with $||\phi''||_{2,\al}\leq c_0\,\e^{3/4}$, 
the problem (\ref{eq:8.4}) has a unique solution $w$. The mapping 
\[
\Pi''\left( {\cal C}^{2,\al}(S^1)\right)^k \ni \phi'' \longrightarrow w \in 
{\cal C}^{2,\al}_{\mu}(\tilde{\Sigma}_\e),
\]
is continuous and the solution $w$ satisfies the estimates
\begin{equation}
||w ||_{2,\al,\mu}\leq c\, \e^{\mu/4}( \e^{3/2} + ||\phi''||_{2,\al}  
+ \e^{-3/4}  ||\phi''||^2_{2,\al})
\label{tdelrefest1}
\end{equation}
and
\[
|| (w - \Pi''w) (S_{\e_\ell}/8,\cdot)||_{2,\al}+ ||\del_s w - \Pi'' w) ( S_{\e_\ell}/8,\cdot) ||_{1,\al} 
\]
\begin{equation}
\leq c \,  ( \e^{3/2} + (\e^{\kappa- 1} + \e^{(2-\mu)/4}) ||\phi''||_{2,\al}  
+ \e^{-3/4}  ||\phi''||^2_{2,\al}).
\label{tdelrefest2}
\end{equation}
Finally, if $\tilde{w}_0 = \tilde{P}_{0}(\phi'') \in {\mathcal C}^{2,\al}_{\mu
}(\tilde{\Sigma}_\e)$  as in Corollary~\ref{cor:8.2}, then 
\begin{equation}
||w - \tilde{w}_0 ||_{2,\al, \mu} \leq c\, \e^{-\mu/4}
( \e^{3/2} + (\e^{\kappa- 1} 
+ \e^{(2-\mu)/4}) ||\phi''||_{2,\al}  + \e^{-3/4}  ||\phi''||^2_{2,\al}).
\label{tdelrefest3}
\end{equation}
\label{cor:2222}
\end{corollary}
The proof is  identical to that of Corollary \ref{cor:1111}, and so 
we omit it.

\section{Matching the Cauchy data} 

We have now established that, given any set of parameters
${\cal P} = (\tilde{t},\tilde{r},\tilde{d},\tilde{\delta})$
satisfying $\|{\cal P}\| \leq \e^\kappa$, and for any 
$\phi'' \in \Pi''\left( {\cal C}^{2,\al}(S^1)\right)^k$, we can
solve the equations
\begin{equation}
\left\{ \begin{array}{rlll}
L_{\e_\ell +\delta^\ell} w_\ell & = & {{\cal Q}} (w_\ell)    
& \qquad \mbox{in} \qquad [S_\e /8, +\infty)\times S^1 \\[2mm]
\Pi'' w_\ell  & = &  \phi''_\ell & \qquad \mbox{on} \qquad 
\{ S_\e /8\} \times S^1 
        \end{array}
\right.
\label{eq:9.2}
\end{equation}
and 
\begin{equation}
\left\{ \begin{array}{rlll}
\tLeP w_{\cal K} & = & 1 + \tilde{{\cal Q}} (w_0) & 
\qquad \mbox{in}  \qquad \tSeP\\[2mm]
\Pi'' w_{\cal K} & = &  \phi'' & \qquad \mbox{on} 
\qquad \del \tSeP,
\end{array}
\right.
\label{eq:9.3}
\end{equation}
when $\e$ is sufficiently small. Thus we may define the mappings
\[
\begin{array}{rlll}
{\cal S}_{\e} :  & \Pi''\left( {\cal C}^{2, \alpha}(S^1) \right)^k 
\ni\phi'' \longrightarrow \\[2mm]
\, & (\del_s w_1(S_{\e_\ell} /8 ,\cdot) ,\ldots ,\del_s w_k( S_{\e_\ell} 
/8 ,\cdot))   \in \left( {\cal C}^{1, \alpha}(S^1) \right)^k ,
\end{array}
\]
and 
\[
\begin{array}{rlll}
{\cal T}_{\e,{\cal P}} : & \Pi'' \left( {\cal C}^{2, \alpha}(S^1)  \right)^k 
\ni \phi''  \longrightarrow \\[2mm]
\, & (\del_s (\hat{w}_0 + w_{\cal K})|_{\e E_1}(S_{\e_\ell} /8 ,\cdot), 
\ldots , \del_s (\hat{w}_0 + w_{\cal K})|_{\e E_k}(S_{\e_\ell} /8 ,\cdot))  
\in \left( {\cal C}^{1, \alpha}(S^1) \right)^k ,
\end{array}
\]
where $\hat{w}_0$ is defined in Proposition \ref{pr:modif}.
These would be the Dirichlet-to-Neumann mappings for the
two nonlinear problems (\ref{eq:9.2}) and (\ref{eq:9.3}), save
for the fact that the low eigencomponents of the Dirichlet
data are not specified. 

It follows from Corollary \ref{cor:1111} and Corollary \ref{cor:2222} that,
for $c_0 >0$ small enough, these mappings are well defined 
from the ball of radius $c_0 \, \e^{3/4}$ in $\left({\cal C}^{2,\al}(S^1)
\right)^k$ into the ball of radius $c \, c_0  \, \e^{3/4}$ in the space 
$\left({\cal C}^{1, \al}(S^1)\right)^k$.  
\begin{proposition}
There exists a constant $\tilde{c}_0 >0$ such that in the ball of radius 
$\tilde{c}_0 \, \e^{3/4}$ in $\left({\cal C}^{2,\al}(S^1)\right)^k$ 
there is a unique $\phi''$  which satisfies the equation 
\begin{equation}
\Pi'' \left[ {\cal S}_{\e,{\cal P}} (\phi'' )\right] =
\Pi''\left[{\cal T}_{\e,{\cal P}} (\phi'') \right]
\label{eq:9.4}
\end{equation}
This solution satisfies
\[
\|\phi''\|_{2,\al} \leq c \, (\e^{2\kappa-1} + \e^{3/2}) ,
\]
\end{proposition}
for some constant $c > 0$ independent of $\e$. 

{\bf Proof :} By Corollaries \ref{cor:1111},  \ref{cor:2222} 
and Proposition \ref{pr:modif}, we may approximate these
partial Dirichlet-to-Neumann maps by the corresponding
maps for the Laplacian on the cylinder $\RR \times S^1$.
More specifically, 
\[
\begin{array}{rlll}
\Pi'' \left[ {\cal S}_{\e,{\cal P}}(\phi'' )\right] & = &
(\del_s w_0^1 (S_\e /8, \cdot) ,\ldots ,\del_s w^k_0 (S_{\e_\ell} /8 ,\cdot))
 \\[2mm] 
& + & O( \e^{1/2} +\e^{(6-3\mu)/4}) ||\phi''||_{2, \alpha} 
+ O(\e^{-3/4}) ||\phi''||^2_{2,\alpha},
\end{array}
\]
and similarly
\[
\begin{array}{rlll}
\Pi'' \left[ {\cal T}_{\e,{\cal P}}(\phi'')\right] & = &
(\del_s \tilde{w}_0|_{\e E_1} (S_\e /8, \cdot) ,\ldots ,\del_s 
\tilde{w}_0 |_{\e E_k} (S_{\e_\ell} /8 ,\cdot))  \\[2mm]
& + & O(\e^{3/2} +\e^{2\kappa -1}) + O( \e^{\kappa - 1} +\e^{(2-\mu)/4}) 
||\phi''||_{2,\al} + O(\e^{-3/4})||\phi''||^2_{2,\al},
\end{array}
\]
where $w^\ell_0 =P_{\e_\ell}(\phi''_\ell)$ and 
$\tilde{w}_0 = P_{\e,{\cal P}}(\phi'')$.  We are using 
here that $\Pi'' \hat{w}_0 = O(\e^{3/2} +\e^{2\kappa -1})$.
Therefore we must find $\phi''$ such that
\[
(\del_s (w_0^1 - \tilde{w}_0|_{\e E_1}) (S_{\e_\ell} /8, \cdot) ,\ldots ,
 \del_s (w^k_0 - \tilde{w}_0|_{\e E_k}) (S_{\e_\ell} /8 ,\cdot)= F(\e,
{\cal P},\phi''),
\]
where 
\[
F = O(\e^{3/2} +\e^{2\kappa -1}) + O(\e^{1/2} +\e^{\kappa -1} + 
\e^{(2-\mu)/4})||\phi''||_{2,\al} + O(\e^{-3/4}) ||\phi''||^2_{2,\al}.
\]

To see that this equation has a solution, we first note that
the corresponding homogeneous linear problem has a unique
solution within the range of $\Pi''$, namely $\phi'' = 0$. 
The linear subspaces of Cauchy data for the Laplacian
on the two half-cylinders $(-\infty,S_\e/8] \times S^1$ and
$[S_\e/8,\infty) \times S^1$, restricted to the range
of $\Pi''$, form a transversal Fredholm pair. This implies
that 
\[
S_0 - T_{0,{\cal P}}
\]
is invertible, where $S_0$ and $T_{0,{\cal P}}$ are the 
Dirichlet-to-Neumann operators for these linear problems.
It is also true that this operator is a first order
pseudodifferential operator. This implies that its
inverse is compact. Our problem now reduces to
\[
\phi'' = (S_0 - T_{0,{\cal P}})^{-1}F(\e,{\cal P}, \phi'').
\]
The operator on the right is compact, by the remarks above.
Furthermore, for $\e$ suficiently small, it maps the
ball of radius $c_0 \e^{3/4}$ to itself, because all the
terms in $F$ decay faster (in $\e$) than $\e^{3/4}$.
Hence, for $\e$ sufficiently small, for every choice
of parameter set ${\cal P}$ with $\|{\cal P}\| \leq \e^\kappa$,
this map must have a unique fixed point. \hfill $\Box$

This proposition allows us to reduce the problem, at last,
to a finite dimensional one. For every $\e$ sufficiently
small and ${\cal P}$ with $\|{\cal P}\| \leq \e^\kappa$,
associate the unique element $\phi''$ in the range of
$\Pi''({\cal C}^{2,\al}(S^1))$. The equation
\[
{\cal S}_\e (\phi'') = {\cal T}_{\e,{\cal P}}(\phi'')
\]
reduces to a system of $k$ nonlinear equations of the form
\[
( -\displaystyle \frac{1}{\cosh S_{\e_\ell}} ( t_1^\ell 
\cos \theta + t_2^\ell \sin \theta )  -  ( r_1^{\ell} \cos \theta + 
r^{\ell}_2 \sin \theta ) \e_\ell  \cosh  S_{\e_\ell}   + d^{\ell}  
 +\delta^\ell  S_{\e_\ell},
\]
\[
\displaystyle \frac{\sinh S_{\e_\ell} }{\cosh^2 S_{\e_\ell}} ( t_1^\ell 
\cos \theta + t_2^\ell \sin \theta )  -  ( r_1^{\ell} \cos \theta + 
r^{\ell}_2 \sin \theta ) \e_\ell  \sinh  S_{\e_\ell} +\delta^\ell)
\]
\[
= O(\e^{3/2} + \e^{2\kappa -1}).
\]
Because of the restriction on the norm of ${\cal P}$, and because
\[
\e^{3/2} + \e^{2\kappa -1} =o(\e^{\kappa}),
\]
we may again conclude that this equation has a solution. 
This ends the proof of our result.

\section{The nondegeneracy of the solutions}

We now show that for $\e$ sufficiently small, the solutions we have 
constructed above are nondegenerate in the sense defined in \cite{KMP}. 
This condition ensures the smoothness of the moduli spaces 
${\cal M}_{g,k}$ near $\Sigma_\e$. We begin by recalling this 
notion of nondegeneracy. 

\begin{definition} The constant mean curvature surface $\Sigma_\e \in 
{\cal M}_{g,k}$ is {\it nondegenerate} if the linearization of the mean 
curvature operator  about $\Sigma_\e$ is injective on the function space 
${\cal C}^{2, \al}_{\delta}(\Sigma_\e)$ for all $\delta < 0 $.
\end{definition}

Here for $r\in {\NN}$, $\al \in [0,1)$ and  $\delta \in {\RR}$, 
${\cal C}^{r, \al}_{\delta} (\Sigma_\e)$ is defined to be the space 
of functions $\phi \in {\cal C}^{r, \al}(\Sigma_\e)$ which can be 
written on each end of $\Sigma_\e$ as $e^{\delta s}$ times a function 
$\psi$ with $\psi \in {\cal C}^{r, \al} ({\RR}^+_s \times S^1_\theta)$.

First notice that it is sufficient to prove that, for $\e$ small enough, 
the Jacobi operator $L$ is injective on ${\cal C}^{2,\al}_{\delta}
(\Sigma_\e)$ for some fixed $\delta \in (-2,-1)$. This is because
any decaying solution of $Lu = 0$ must decay exponentially near the 
$i$th end of $\Sigma_\e$ at least like $e^{-\gamma_2(\e_i)s}$, and
by Proposition~\ref{pr:4.15}, when $\e$ is sufficiently small,
$2 - \gamma_2(\e_i)$ is as small as desired, so that 
$u \in {\cal C}^{2, \al}_{\delta}(\Sigma_\e)$.

The proof is by contradiction. Fix $\delta \in (-2,-1)$ and assume that 
for some sequence of $\e_k$ tending to $0$, the Jacobi operator 
\[
{\cal L}_k = \Delta_{\Sigma_{\e_k}} + |{\bf A}_{\Sigma_{\e_k}}|^2
\]
on $\Sigma_{\e_k}$ is not injective on ${\cal C}^{2,\al}_{\delta}
(\Sigma_{\e_k})$. Then there exists some $w_k \in {\cal C}^{2,\al}_{\delta}
(\Sigma_{\e_k})$ such that ${\cal L}_k w_k = 0$ and $w_k \neq 0$.

First normalize $w_k$, multiplying it by a suitable constant, so that
$||w_k ||_{0,0,\delta} (\Sigma_{\e_k}) = 1$.  Choose a point $y_k \in 
\Sigma_{\e_k}$ where the above norm is achieved. Suppose first that some 
subsequence of the $y_k/ \e_k$ converges to a point $y_0 \in \Sigma_0$. 
Then we can extract a subsequence of the $w_k$ which converge on every 
compact of $\Sigma_0$ to a limiting function $w$ globally defined 
on $\Sigma_0$; $w$ must be nontrivial since we also have $||w||_{0,0, 
\delta (\Sigma_0)}=1$. 
Furthermore, $L_{\Sigma_0} w  = 0$. Since we have asssumed that $\Sigma_0$ is 
nondegenerate, we have obtained a contradiction.

If, on the other hand, some subsequence of the $y_k$ satisfies 
$\lim_{k\rightarrow +\infty} |y_k/\e_k|=+\infty$ then, this implies 
that, at least for a subsequence, the points $y_k$ are always in the
 same end, say the $i$th. Therefore, we may write, 
\[
y_k = \x_{\e_k,i} (s_k + s_{\e_k,i} , \theta_k),
\]
with $s_k$ tending to $+\infty$. By translating back by $s_k + s_{\e_k,i}$ and
multiplying by a suitable constant, we find yet another sequence
of solutions, which we again call $w_k$, attaining their maximum
at $s = 0$, and which solve the translated equation, which we
again write as ${\cal L}_k w_k = 0$. Here ${\cal L}_k$ is the
linearized mean curvature operator relative to the parameterization 
given above near the ends.
It is straightforward to see that the $w_k$ converge 
to a nontrivial solution $w$ of one of the following two limiting equations
\begin{equation}
\del_s^2 w +  \del_\theta^2 w  = 0,
\end{equation}
or  
\begin{equation}
\del_s^2  + \del_\theta^2 w +\frac{2}{\cosh^2 (s+\bar{s})}w= 0,
\qquad \mbox{ for some} \quad\bar{s}\in {\RR} ,
\end{equation}
on ${\RR} \times S^{1}$. In addition, $w$ is bounded by $e^{\delta s}$.
By decomposing into eigenfrequencies we then see that necessarly $w=0$ 
which is the desired contradiction.

This covers all cases, so we have showed that the linearization is
injective on the appropriate weighted H\"older spaces.

 \end{document}